\documentclass[11pt]{article}
\usepackage{authblk}
\usepackage[utf8]{inputenc}
\usepackage{amssymb}
\usepackage{amsmath}
\usepackage{amsthm}
\usepackage{multicol}
\usepackage{graphics}
\usepackage{alltt,listings,enumerate,hhline}
\newtheorem{lemma}{Lemma}
\newtheorem{proposition}[lemma]{Proposition}
\newtheorem{example}[lemma]{Example}
\newtheorem{theorem}[lemma]{Theorem}
\newtheorem{remark}[lemma]{Remark}
\newtheorem{corollary}[lemma]{Corollary}
\usepackage[section]{placeins}
\newcommand\mut[1]{\ignorespaces}
\DeclareMathOperator{\rank}{rank}
\DeclareMathOperator{\intmo}{int}
\usepackage{algorithm,algpseudocode}
\usepackage{tikz,tikz-qtree,tikz-qtree-compat,adjustbox}
\usepackage{pst-plot}
\usepackage{etoolbox}

\algnewcommand{\algorithmicand}{\textbf{ and }}
\algnewcommand{\algorithmicor}{\textbf{ or }}
\algnewcommand{\Or}{\algorithmicor}
\algnewcommand{\And}{\algorithmicand}
\algnewcommand{\Var}{\texttt}

\renewcommand\leq\leqslant
\renewcommand\geq\geqslant

\usepackage{tikz,tikz-qtree,tikz-qtree-compat,fullpage,colortbl}

\title{On the seeds and the great-grandchildren\\ of a numerical semigroup}
\date{\today}

\author[$\star$]{Maria Bras-Amorós}

\affil[$\star$]{Universitat Rovira i Virgili, Av. Països Catalans 26, 43007, Tarragona, Catalonia, Spain. \small\textit{E-mail adress:} \texttt{maria.bras@urv.cat}}

\begin{document}
\maketitle

\begin{abstract}
  We present a revisit of the seeds algorithm to explore the semigroup tree.
  First, an equivalent definition of seed is presented, which seems easier to manage. Second, we determine the seeds of semigroups with at most three left elements. And third, we find the great-grandchildren of any numerical semigroup in terms of its seeds.
  
  The RGD algorithm is the fastest known algorithm at the moment.
  But if one compares the originary seeds algorithm with the RGD algorithm, one observes that the seeds algorithm uses more elaborated mathematical tools while the RGD algorithm uses data structures that are better adapted to the final C implementations. For genera up to around one half of the maximum size of native integers, the newly defined seeds algorithm performs significantly better than the RGD algorithm. For future compilators allowing larger native sized integers this may constitute a powerful tool to explore the semigroup tree up to genera never explored before. The new seeds algorithm uses bitwise integer operations, the knowledge of the seeds of semigroups with at most three left elements and of the great-grandchildren of any numerical semigroup, apart from techniques such as parallelization and depth first search as wisely introduced in this context by Fromentin and Hivert.
 
  The algorithm has been used to prove that there are no Eliahou semigroups of genus $66$, hence proving the Wilf conjecture for genus up to $66$. We also found three Eliahou semigroups of genus $67$. One of these semigroups is neither of Eliahou-Fromentin type, nor of Delgado's type. However, it is a member of a new family suggested by Shalom Eliahou.
\end{abstract}


\section{Introduction}
A {\em numerical semigroup} is a cofinite submonoid of ${\mathbb N}$.
The elements in its complement in ${\mathbb N}$ are denoted the {\em gaps}
of the semigroup and the {\em genus} of the semigroup is the number of its gaps. There is a finite number of semigroups of each given genus. A major problem in the theory of numerical semigroups is the exhaustive listing of all numerical semigroups for increasing values of the genus. See \cite{RGS} for a general reference on numerical semigroups.

The {\em primitive elements} (or minimal generators) of a numerical semigroup are those elements of the semigroup that can not be obtained as a sum of two smaller semigroup elements. If we take away a primitive element from a numerical semigroup we obtain another semigroup with genus increased by one.
The elements of the semigroup that are smaller than some gap are denoted the {\em left elements} of the semigroup.
The primitive elements that are larger than the largest gap are called {\em right generators}.
We can organize all numerical semigroups in an infinite tree rooted at ${\mathbb N}$ and such that the children of a node are the semigroups obtained taking away one by one its right generators.

The tree of semigroups was first extensively explored in \cite{Br:fibonacci}, in this case using brute force to find right generators. In \cite{Br:bounds} it was noticed that the search of right generators of a child can be restricted to the right generators of its parent and to one further element. Indeed, suppose that $\lambda_1$ is the smallest nonzero non-gap (i.e. the {\em multiplicity}) of a numerical semigroup $\Lambda$. If the gaps of $\Lambda$ are not all in a row (i.e. if $\Lambda$ is not {\em ordinary}), when we take away one right generator $\lambda$ of $\Lambda$, the right generators of the new semigroup $\Lambda\setminus\{\lambda\}$ all belong to the set of right generators of $\Lambda$ except the element $\lambda+\lambda_1$, should it be a primitive element of the new semigroup. If $\lambda+\lambda_1$ is a primitive element of the child, then $\lambda$ is called a {\em strong} generator of $\Lambda$. Otherwise it is called a {\em weak} generator of $\Lambda$. Fromentin and Hivert presented in \cite{FromentinHivert,FH-code} a very efficient algorithm using parallel computation and depth first search, representing semigroups by their decomposition numbers. Suppose that the numerical semigroup $\Lambda$ is $\{\lambda_0=0<\lambda_1<\dots\}$. The sequence of decomposition numbers of $\Lambda$ can be defined by $d_i=\left\lceil\frac{\#\{\lambda_j\in\Lambda:\lambda_i-\lambda_j\in\Lambda\}}{2}\right\rceil$. See its relationship with the $\nu$ sequence in \cite[Exercise 2.16]{RGS} and references therein.
In \cite{seeds} an algorithm was presented using the notion of {\em seed}.
An element $\lambda_s$ larger than the largest gap is a $p$-th order seed of $\Lambda$ if $\lambda_s+\lambda_p\neq \lambda_i+\lambda_j$ for all $p<i\leq j<s$.
Now, right generators are order-$0$ seeds and they are strong generators if and only if they are order-$1$ seeds.
From another perspective, in \cite{rgd,RGD-code} a very simple but very efficient algorithm was presented, using the idea of right-generators descendant of a numerical semigroup. The {\em right-generators descendant} of a numerical semigroup $\Lambda$ is the semigroup obtained by taking away from $\Lambda$ {\em all} its right generators. The so-called RGD algorithm is at the moment the most efficient algorithm.

We present a revisit of the seeds algorithm to explore the semigroup tree. First, an equivalent definition of seed is presented, which seems easier to manage. Second, we determine the seeds of semigroups with at most three left elements. And third, we find the great-grandchildren of any numerical semigroup in terms of its seeds.
For genera around one half of the maximum size of native integers, the newly defined seeds algorithm performs significantly better than the RGD algorithm. For future compilators allowing larger native sized integers this may constitute a powerful tool to explore the semigroup tree up to genera never explored before. The new seeds algorithm uses bitwise integer operations apart from techniques such as parallelization and depth first search as wisely introduced in this context by Fromentin and Hivert.

Using 128-bit integers, as is the current bit length limitation in C, one can explore semigroups up to genus $64$, search for Eliahou semigroups and check the Wilf conjecture for genus up to $66$ and count semigroups of genus up to $69$. In an eventual $2b$-bit integer limitation scenario, our algorithm will be able to efficiently explore semigroups up to genus $b$, search for Eliahou semigroups and check the Wilf conjecture for genus up to $b+2$ and count semigroups of genus up to $b+5$. In particular, the algorithm has been used to prove that there are no Eliahou semigroups of genus $66$, hence proving the Wilf conjecture for genus up to $66$.

In Section~\ref{s:bitstreams} we revisit the definition of the gap bitstream and the seed bitstream of a numerical semigroup, this time using an alternative approach to the definition of seed with respect to the one in \cite{seeds}. In
Section~\ref{s:smallrank} we determine these bitstreams for semigroups of small rank. In Section~\ref{s:bitstreamsupdate} we revisit how the seed bistream of a semigroup can be computed from the seed bitstrem of its parent from the new perspective. In Section~\ref{s:specialcases} we characterize new generators, new strong generators and newly strong generators from an analysis of the seed bitsream. In Section~\ref{s:greatgrandchildren} we explain how the bitstreams can be used to determine the great-grandchildren of a numerical semigroup. All the previous results can be used to build a revisited version of the seeds algorithm and in Section~\ref{s:performance} we present empiric performance results. In Section~\ref{s:eliahou} we recall the Wilf conjecture, and Eliahou and Delgado's related results and use our algorithm to check the Wilf conjecture for genus 66 and to find three new Eliahou semigroups of genus $67$, one is of Eliahou-Fromentin type, another one is of Delgado's type and the third one is neither of Eliahou-Fromentin type, nor of Delgado's type.
However, we will show that it is a member of a new family suggested by Shalom Eliahou.
In Section~\ref{s:howfar} we analyze the limits of our algorithm in terms of the integer bit size bound of the compilator in which the algorithm is implemented. In particular it is analyzed the limit of our experimental results in a C compilator with a bound of 128 bits per integer.

\section{The bitstream of gaps and the bitstream of seeds of a numerical semigroup}\label{s:bitstreams}

If a numerical semigroup $\Lambda$ is $\{\lambda_0=0<\lambda_1<\dots\}$, define its {\em multiplicity} as $m(\Lambda)=\lambda_1$, its Frobenius number, $F(\Lambda)$, as its largest gap and its {\em conductor}, $c(\Lambda)$, as its largest gap plus one.
Let $L(\Lambda)$ be the set of left elements of $\Lambda$.
A numerical semigroup is called {\em ordinary} if $c(\Lambda)=m(\Lambda)=\lambda_1$ and 
{\em pseudo-ordinary} if $c(\Lambda)=\lambda_2$.
We say that the {\em rank} of $\Lambda$ is $$k(\Lambda):=c(\Lambda)-g(\Lambda).$$
This way, $\lambda_{k(\Lambda)}=c(\Lambda)$ and so the rank of $\Lambda$ coincides with the cardinal of $L(\Lambda)$.
Now, the unique semigroup of rank $0$ is ${\mathbb N}$, 
the semigroups of rank $1$ are the non-trivial ordinary semigroups, and the semigroups of rank $2$ are pseudo-ordinary semigroups. 
Define the {\em jump function} as
$$\begin{array}{rcccl}
  u: & {\mathbb N} & \longrightarrow & {\mathbb N} &\\
     & j           & \mapsto &u_j&:=\lambda_{j+1}-\lambda_j\\
  \end{array}$$
Observe that $u_j=1$ for all $j\geq k(\Lambda)$, and that $\sum_{j=0}^{k(\Lambda)-1}u_j=c(\Lambda)$. Furthermore, the multiplicity coincides with $u_0$ and $u_j\leq u_0$ for all $j\in{\mathbb N}$.

A bitstream is a finite sequence $a=a_0\dots a_\ell$ where $a_i$ is either $0$ or $1$ for every $i$.
For our purposes, we can indistinctly use $a$ for $a_0\dots a_\ell$ and for any bitstream of the form $a_0\dots a_\ell\underbrace{0\dots0}_{t}$ for any positive integer $t$.
The {\em weight} of a bitstream is the number of its nonzero elements. The weight of a bitstream $a$ will be denoted $w(a)$. The weight of the substream $a_i\dots a_j$ of $a$ will be denoted $w_i^j(a)$.
If $n\in{\mathbb N}$ we denote $(n)_b$ its binary representation.
That is, $(n)_b=a_0\dots a_\ell$ where $a_0,\dots,a_\ell$ are the unique values in $\{0,1\}$ such that $n=a_0+2a_1+2^2a_2+\dots+2^\ell a_\ell$.
For instance, $$(2^i)_b=\underbrace{0\dots0}_{i}1,\ \ \ \ \ \ \ \ \ \ \ \ \ \ \ \ (2^i-1)_b=\underbrace{1\dots1}_{i}.$$ On the contray, for a bitstream $a=a_0,\dots,a_\ell$ we denote $\intmo(a)$ the integer $a_0+2a_1+2^2a_2+\dots+2^\ell a_\ell$.

Suppose a semigroup $\Lambda$ has conductor $c$. From now on, whenever the context is clear, we will omit $\Lambda$ in $c(\Lambda),k(\Lambda)\dots$. We encode its gaps as the bitstream $$G(\Lambda)=g_0\dots g_{c-1}$$ with $g_i=0$ if $i+1\in \Lambda$ and $g_i=1$ otherwise.
Notice that $$G(\Lambda)=\frac{1}{2}\left(2^{c(\Lambda)}-1-\sum_{j=0}^{k(\Lambda)-1}2^{\lambda_j}\right)_b.$$

Recall that for $p<k(\Lambda)$, an element $\lambda_s$ larger than or equal to $c$ is an order-$p$ {\em seed} of $\Lambda$ if $\lambda_s+\lambda_p\neq \lambda_i+\lambda_j$ for all $p<i\leq j<s$. One can check
that any order-$p$ seed is at most $c + u_p- 1$. In particular,
the number of order-$p$ seeds is at most $u_p$ 
(\cite[Lemma 1.2.]{seeds}).

In \cite{seeds}, the {\em table of seeds} of a semigroup $\Lambda$ was defined as the binary table whose rows are indexed in $0 \leq i \leq k(\Lambda) - 1$, the $i$th row has $u_i$ entries, and the $j$th entry in the $i$th row (with $0\leq j<u_i$) is $1$ if $c + j$ is an order-$i$ seed of $\Lambda$, and $0$ otherwise.

\begin{example}
  The table of seeds of the semigroup $\Lambda=\{0,3,6,8,9,10,\dots\}$, with conductor $c=8$, is
  $$\begin{array}{|c|c|c|}
    \hline
    1&0&1\\\hline
    1&0&1\\\hline
    1&1&\multicolumn{1}{c}{}\\\cline{1-2}
    \end{array}$$
\end{example}

We can also encode the seeds of $\Lambda$ as the bitstream $$S(\Lambda)=S_0\dots S_{c-1}$$with $S_i=1$ if $c+i-\lambda_j$ is a $j$-th order seed of $\Lambda$ where $j$ is the unique non-negative integer such that $\lambda_j\leq i <\lambda_{j+1}$.
In fact, the seed bitstream is nothing else but the concatenation of the rows of the table of seeds.

\begin{example}
  The same semigroup $\Lambda=\{0,3,6,8,9,10,\dots\}$ has $G(\Lambda)=11011010$ and $S(\Lambda)=10110111$.
\end{example}

\begin{figure}
\input{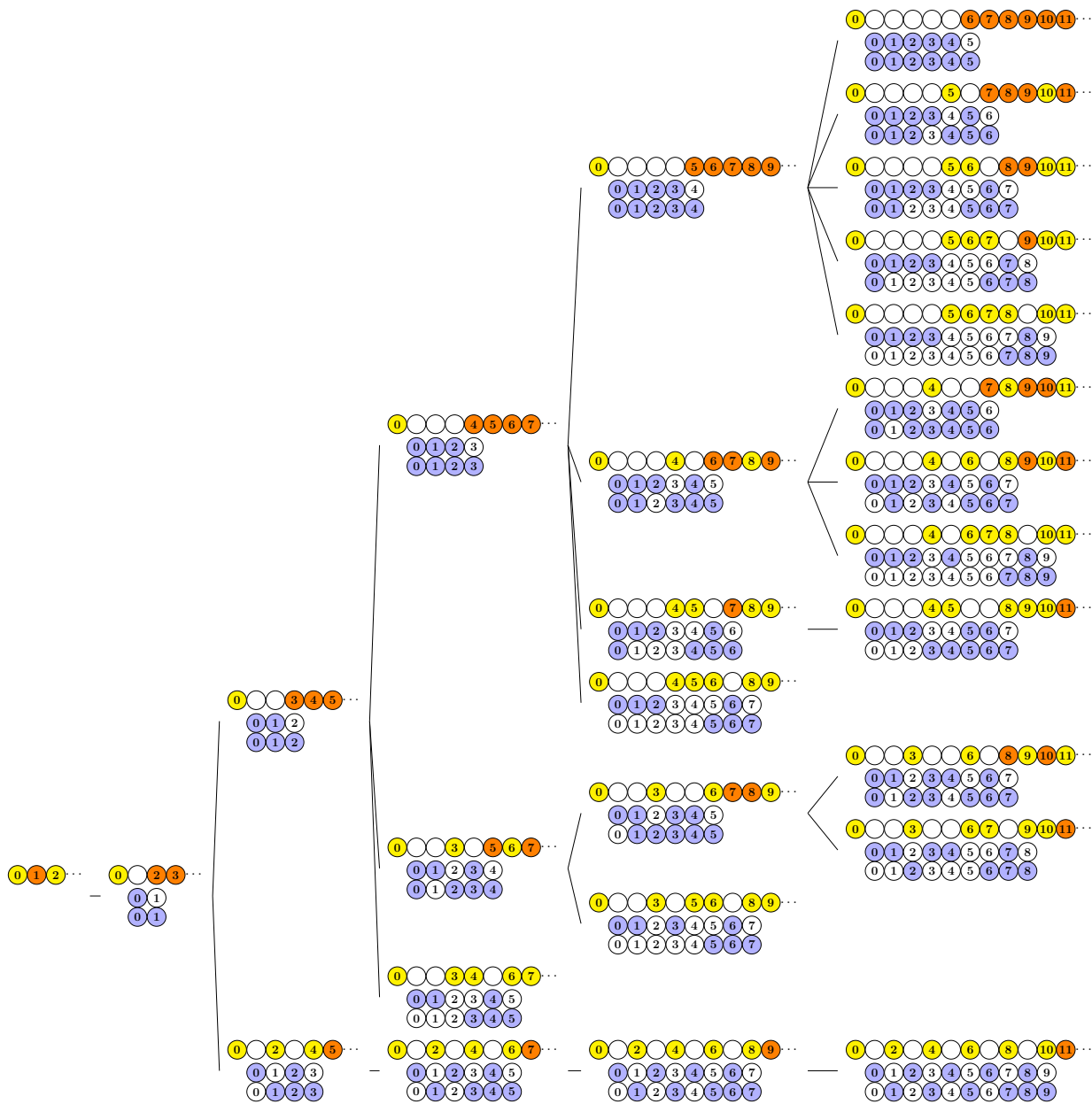}
\caption{Lowest depth nodes of the semigroup tree. Each node is represented by its gap and its seed sequences.}
\label{f:tree}
\end{figure}

In Figure~\ref{f:tree} one can see the sequences $G$ and $S$ for the lowest genus semigroups organized in the semigroup tree. To make it easier to read, we represented each $1$ in the sequences with a dark circle with its position written inside, and each $0$ in the sequence with a light gray circle with its position also written inside.

\begin{lemma}
  Given an integer $\lambda$ with $c\leq \lambda< 2c$, let $j\in\{0,\dots,k-1\}$ be such that $\lambda_j\leq \lambda-c<\lambda_{j+1}$. Then,
  \begin{enumerate}[(a)]
\item
  $\lambda-\lambda_j$ is an order-$j$ seed of $\Lambda$ if and only if
  $\lambda$ is not the sum of two elements in $L(\Lambda)$.
  \item   $\lambda-\lambda_j$ is an order-$j$ seed of $\Lambda$ if and only if
  $\lambda$ is not the sum of two elements in $\{\lambda_{j+1},\dots,\lambda_{k-1}\}$.
  \end{enumerate}
  \end{lemma}

\begin{proof}
Note that $\lambda=\lambda_s$ with $s=\lambda-c+k(\Lambda)$.
  The non-gap $\lambda-\lambda_j$ is {\em not} an order-$j$ seed of $\Lambda$ if and only if
  $\lambda=(\lambda-\lambda_j)+\lambda_j$ is the sum of two elements in $\{\lambda_{j+1},\dots,\lambda_{s-1}\}$.
  As $\lambda$ can not be the sum $\lambda_u+\lambda_v$ with $\lambda_u>0$ and $v\geq s$, we obtain that
   $\lambda-\lambda_j$ is {\em not} an order-$j$ seed of $\Lambda$ if and only if $\lambda$ is the sum of two elements in $\{\lambda_{j+1},\lambda_{j+2},\dots\}$.
In this case $$\lambda=\lambda_{j+1+\ell_1}+\lambda_{j+1+\ell_2}$$ for some non-negative integers $\ell_1$ and $\ell_2$. Each of the two summands can now be bounded as follows.
  \begin{itemize}
  \item $\lambda_{j+1+\ell_1}<c$ since
    $\lambda_{j+1+\ell_1}=\lambda-\lambda_{j+1+\ell_2}\leq \lambda-\lambda_{j+1}<c$.
  \item $\lambda_{j+1+\ell_2}<c$ since
    $\lambda_{j+1+\ell_2}=\lambda-\lambda_{j+1+\ell_1}\leq \lambda-\lambda_{j+1}<c$.
  \end{itemize}
\end{proof}

Next corollary gives an alternative definition of the seed bitstream, and hence, an alternative definition of the notion of seed.

\begin{corollary}\label{cor:sigma}
  Let $\Sigma(\Lambda)$ be the bitstream $\Sigma_0\Sigma_1\dots\Sigma_{2c-1}$, with
$\Sigma_i=0$ if $i$ belongs to $L(\Lambda)+L(\Lambda)$ and $\Sigma_i=1$ otherwise,
then
$$G_0G_1\dots G_{c-2}=\Sigma_{1}\dots\Sigma_{c-1},\ \ \ \ \ G_{c-1}=0$$
$$S_0S_1\dots S_{c-1}=\Sigma_{c}\Sigma_{c+1}\dots\Sigma_{2c-1},$$
Hence, $$G=\left((\intmo(\Sigma)\rm{\ mod\ }2^c)//2\right)_b,$$ $$S=\left(\intmo(\Sigma)//2^{c}\right)_b,$$
where $a//b$ and $a\mod b$ denote the quotient and the remainder of the division of $a$ by $b$, respectively.
\end{corollary}

\begin{example}
  Following the same example $\Lambda=\{0,3,6,8,9,10,\dots\}$, we have
$L(\Lambda)=\{0,3,6\}$ and $L(\Lambda)+L(\Lambda)=\{0,3,6,9,12\}$, that is $\Sigma=0110110110110111$.
$$\overbrace{0\underbrace{1\,1\,0\,1\,1\,0\,1}_{G(\Lambda)_{0,\dots,(c-2)}}\underbrace{1\,0\,1\,1\,0\,1\,1\,1}_{S(\Lambda)}}^{\Sigma(\Lambda)}$$
\end{example}

As a consequence of this alternative definition of $S(\Lambda)$ together with the fact that 
$\max(L(\Lambda)+L(\Lambda))=2\lambda_{k-1}=2c-2u_{k-1}$ we derive that 
$\Sigma_{2\lambda_{k-1}}=0$ and
$\Sigma_{2\lambda_{k-1}+1}=\Sigma_{2\lambda_{k-1}+2}=\dots=\Sigma_{2c-1}=1$. That is,
$$S=S_0\dots S_{2\lambda_{k-1}-1}\ \ \ 0\underbrace{1\cdots1}_{2u_{k-1}-1}$$

\begin{corollary}
  The elements $c, c+1, \dots, c+u_{k-1}-1$ are order-$(k-1)$ seeds of $\Lambda$.
\end{corollary}

\mut{
\begin{remark}
Although this will not be used in the algorithm exposed at the end, the seed bitstream can be obtained from the gap bitstream as follows.

\begin{alltt}
  for(i=0;i<c;i++)
  \{
    j=i;
    while(j>0 && \(G\sb{j-1}\)==1)
      j--;
    \(S\sb{i}\)=1;
    for(k=1; k<=(c+i-2*j)/2 && \(S\sb{i}\)==1; k++)
      \{
         if(\(G\sb{i-j+c-k-1}\)=0 && \(G\sb{j+k-1}\)==0)
           \(S\sb{i}\)=0;
      \}
  \}
\end{alltt}
\end{remark}
}

\section{The gap and the seed bitstreams of semigroups with small rank}
\label{s:smallrank}

Next lemma states the gap bitstream and the seed bitstream for semigroups of rank $k$ at most $3$.
We will call $u:=u_1$ and $v:=u_2$.


\begin{lemma}\label{l:lowrank}
The gap and the seed bitsreams of semigroups of rank $k$ between $1$ and $3$ are given by the following formulas.  
  \begin{itemize}
  \item For semigroups of rank $k=1$, 
    \begin{eqnarray*}G(\{0,m,\dots\})&=&(2^{m-1}-1)_b\\S(\{0,m,\dots\})&=&(2^m-1)_b\end{eqnarray*}
  \item For semigroups of rank $k=2$ (in which case $1<u\leq m$), 
    \begin{eqnarray*}G(\{0,m,m+u,\dots\})&=&(2^{m+u-1}-1-2^{m-1})_b\\S(\{0,m,m+u,\dots\})&=&(2^{m+u}-1-2^{m-u})_b\end{eqnarray*}
    \item For semigroups of rank $k=3$,
      \begin{itemize}
\item        If $1\leq u<m$,
        then $2\leq v\leq m-u$ and
      \begin{eqnarray*}G(\{0,m,m+u,m+u+v,\dots\})&=&(2^{m+u+v-1}-1-2^{m-1}-2^{m+u-1})_b
        \\S(\{0,m,m+u,m+u+v,\dots\})&=&(2^{m+u+v}-1-2^{m-u-v}-2^{m-v}-2^{m+u-v})_b\end{eqnarray*}        
\item If $u=m$, then $2\leq v\leq m$ and
  \begin{eqnarray*}G(\{0,m,m+u,m+u+v,\dots\})&=&(2^{m+u+v-1}-1-2^{m-1}-2^{m+u-1})_b
    \\S(\{0,m,m+u,m+u+v,\dots\})&=&(2^{m+u+v}-1-2^{m-v}-2^{m+u-v})_b\end{eqnarray*}
  \end{itemize}\end{itemize}
  \end{lemma}

\begin{proof} 
  \begin{itemize}
  \item
    For the rank-$1$ case $L(\Lambda)+L(\Lambda)=\{0\}$ and so, $$\Sigma=(2^{2c}-2)_b=0\underbrace{1\dots 1}_{2c-1}.$$ Then the result follows from Corollary~\ref{cor:sigma}.
  \item
    For the rank-$2$ case, $L(\Lambda)+L(\Lambda)=\{0,m,2m\}$ and so, \begin{eqnarray*}\Sigma&=&(2^{2c}-2-2^m-2^{2m})_b\\&=&(2^c(2^c-2^{2m-c}-1)+(2^c-2^m-2))_b.\end{eqnarray*}
Hence, by Corollary~\ref{cor:sigma}, \begin{eqnarray*}G&=&(2^{c-1}-2^{m-1}-1)_b\\&=&(2^{m+u-1}-2^{m-1}-1)_b,\end{eqnarray*} and \begin{eqnarray*}S&=&(2^c-2^{2m-c}-1)_b\\&=&(2^{m+u}-2^{m-u}-1)_b.\end{eqnarray*}

\item
  For the rank-$3$ case suppose $\Lambda=\{0,m,m+u,m+u+v,\dots\}$ with $v>1$. Then the conductor of $\Lambda$ is $c=m+u+v$.

If $u<m$ then $m+u<2m$ and since $2m\in\Lambda$, we have $2m\geq c$ and so
$v\leq m-u$. Then $L(\Lambda)+L(\Lambda)=\{0,m,m+u,2m,2m+u,2m+2u\}$ with $m+u<c\leq 2m$.
So, \begin{eqnarray*}\Sigma&=&(2^{2c}-2-2^m-2^{m+u}-2^{2m}-2^{2m+u}-2^{2m+2u})_b\\&=&(2^c(2^c-2^{2m+2u-c}-2^{2m+u-c}-2^{2m-c}-1)+(2^{c}-2^{m+u}-2^{m}-2))_b\end{eqnarray*} and, by Corollary~\ref{cor:sigma},
\begin{eqnarray*}G&=&(2^{c-1}-2^{m+u-1}-2^{m-1}-1)_b\\&=&(2^{m+u+v-1}-2^{m+u-1}-2^{m-1}-1)_b\end{eqnarray*} and
\begin{eqnarray*}S&=&(2^c-2^{2m+2u-c}-2^{2m+u-c}-2^{2m-c}-1)_b\\&=&(2^{m+u+v}-2^{m+u-v}-2^{m-v}-2^{m-u-v}-1)_b.\end{eqnarray*}

Otherwise,  $u=m$ and $L(\Lambda)+L(\Lambda)=\{0,m,2m,3m,4m\}$ and so, \begin{eqnarray*}\Sigma&=&(2^{2c}-2-2^m-2^{2m}-2^{3m}-2^{4m})_b\\&=&(2^c(2^c-2^{4m-c}-2^{3m-c}-1)+(2^c-2^{2m}-2^m-2))_b.\end{eqnarray*}
Hence, by Corollary~\ref{cor:sigma}, \begin{eqnarray*}G&=&(2^{c-1}-2^{2m-1}-2^{m-1}-1)_b\\&=&(2^{m+u+v-1}-2^{m+u-1}-2^{m-1}-1)_b,\end{eqnarray*} and \begin{eqnarray*}S&=&(2^c-2^{4m-c}-2^{3m-c}-1)_b\\&=&(2^{m+u+v}-2^{m+u-v}-2^{m-v}-1)_b.\end{eqnarray*}
%
%
\end{itemize}
  \end{proof}

\section{Updating the bitstream of seeds from parents to children}
\label{s:bitstreamsupdate}

Suppose that $\Lambda$ is a semigroup of rank $k$ and conductor $c$ and let $\tilde \Lambda=\Lambda\setminus\{\lambda_s\}$ with $s\geq k$.
We will denote $\tilde \Sigma=\Sigma(\tilde \Lambda)$, $\tilde G=G(\tilde \Lambda)$, $\tilde S=S(\tilde \Lambda)$, $\tilde c=c(\tilde\Lambda)=\lambda_s+1$, $\tilde k$ the rank of $\tilde\Lambda$,  $\tilde u_j$ the $j$th jump of $\tilde \Lambda$ and $\tilde\lambda_j$ the $j$th element of $\tilde\Lambda$.
Let us denote $L=L(\Lambda)$, $\tilde L=L(\tilde\Lambda)$, and $I_s$ the interval $[c,\lambda_s-1]$, which may be empty.
We have $\tilde L=L\cup I_s$. Hence,
$$
\tilde L +\tilde L=\left( L+L \right)\cup\left( L+I_s\right)\cup\left(I_s+I_s\right)$$
with
  $L+I_s=\left\{\lambda_i+j:1\leq i<k\mbox{ and }c\leq j<\lambda_s\right\}$ and
  $I_s+I_s=\{\ell: 2c\leq \ell \leq 2\lambda_s-2\}.$
In particular,

$$\tilde\Sigma_\ell=\left\{\begin{array}{ll}
0 & \mbox{ if } \ell<2c \mbox{ and }\ell\in L+I_s\\
\Sigma_\ell & \mbox{ if } \ell<2c \mbox{ and }\ell\not\in L+I_s\\
0 & \mbox{ if } 2c\leq \ell \leq 2\lambda_s-2\\
1 & \mbox{ if } 2\lambda_s-1\leq \ell \leq 2\lambda_s+1
\end{array}\right.$$

And so,

\begin{eqnarray}\tilde S_\ell&=&\tilde\Sigma_{\ell+\lambda_s+1}\nonumber
  \\&=&
\left\{\begin{array}{ll}
0 & \mbox{ if }
\ell+\lambda_s+1<2c \mbox{ and } \ell+\lambda_s+1\in L+I_s\\
\Sigma_{\ell+\lambda_s+1}=S_{\ell+\lambda_s-c+1} &
\mbox{ if }
\ell+\lambda_s+1<2c \mbox{ and } \ell+\lambda_s+1\not\in L+I_s\\
0 & \mbox{ if } 2c\leq \ell+\lambda_s+1 \leq 2\lambda_s-2\\
1 & \mbox{ if } 2\lambda_s-1\leq \ell+\lambda_s+1 \leq 2\lambda_s+1
\end{array}\right.\label{eq:seeds}
\end{eqnarray}

From this formula we can derive Algorithm~\ref{a:primi} to obtain $\tilde G$ and $\tilde S$ from $G,S,c,k,s$.
An intermediate bitstream variable $auxS$ is introduced which can be defined
so that
$\tilde S_\ell=auxS_{\ell+\lambda_s-c+1}$ except for the case $2\lambda_s-1\leq \ell+\lambda_s+1 \leq 2\lambda_s+1$.
That is
\begin{eqnarray*}
  auxS_\iota&=&\left\{
  \begin{array}{ll}
    S_\iota &\mbox{ if }\iota<c \mbox{ and } \iota+c\not\in L+I_s\\
    0 & \mbox{ otherwise.}
  \end{array}\right.
  \end{eqnarray*}
It can be easily computed since
\begin{eqnarray*}
  auxS_\iota&=&\left\{
  \begin{array}{ll}
    0 &\mbox{ if }\iota+c=\lambda_j+\tilde i \mbox{ for some }\lambda_j\in L(\Lambda)\mbox{ and some }\tilde i\mbox{  with }c\leq \tilde i\leq c+s-k-1\\
    S_\iota&\mbox{ otherwise,}
    \end{array}\right.
  \end{eqnarray*}
or, equivalently,
\begin{eqnarray*}
  auxS_\iota&=&\left\{
  \begin{array}{ll}
    0 &\mbox{ if }\iota=\lambda_j+ i \mbox{ for some }\lambda_j\in \Lambda\mbox{ and some }\tilde i\mbox{  with }0\leq  i\leq s-k-1\\
    S_\iota&\mbox{ otherwise.}
    \end{array}\right.
  \end{eqnarray*}
This is what is done in the for loop of the algorithm.

This algorithm to update $G$ and $S$ is the keystone of the seeds algorithm.
In fact, it may seem more natural to compute $\tilde \Sigma$ from $\Sigma$, but the restriction on the length of integers in implementations makes it more useful to split $\Sigma$ into $G$ and $S$ and $\tilde \Sigma$ into $\tilde G$ and $\tilde S$.
We assume that integers are represented from the least significant bit to the most significant bit. For instance, $(1)_b=100000\dots$.
If $a=a_0\dots a_\ell$
we call $a\gg t$ the bitstream $\underbrace{0\dots 0}_ta_0\dots a_\ell$
and $a\ll t$ the bitstream $a_t\dots a_\ell$.
Notice that $(n)_b\gg t=(2^tn)_b$ and $(n)_b\ll t$ is $(n//2^t)_b$, which is equivalent to discarding the first $t$ positions.
If $a=a_0\dots a_\ell$ and $a'=a'_0\dots a'_{\ell'}$, then
the operations $a\land a'$
and $a\lor a'$
are defined as usual,
\begin{center}
$a\land a'=r_0\dots r_{\max{\{\ell,\ell'\}}}$, where
  $r_i=\left\{\begin{array}{ll}1&\mbox{if }a_i=1\mbox{ and }a'_i=1\\0&\mbox{otherwise}\end{array}\right.$

  $a\lor a'=r_0\dots r_{\max{\{\ell,\ell'\}}}$, where
$r_i=\left\{\begin{array}{ll}1&\mbox{if }a_i=1\mbox{ or }a'_i=1\\0&\mbox{otherwise}\end{array}\right.$
  \end{center}

\begin{algorithm}
  \caption{Obtain $\tilde G$ and $\tilde S$ from $G,S,c,k,s$.}\label{a:primi}
\textbf{Input:} $G=G(\Lambda)$, $S=S(\Lambda)$, $c=c(\Lambda)$, $k=k(\Lambda)$, $s$
\begin{algorithmic}
  \State $A=G$
  \State $auxS=S$
\For{$i$ {\bf from} $0$ {\bf to} $s-k-1$}
\State $A=A\gg 1$
\State $auxS=auxS\land A$
\EndFor
\State$\tilde G=G\lor(1\gg(c+s-1))$
\State$\tilde S=(auxS\ll(s-k+1))\lor (7\gg (c+s-2))$
\end{algorithmic}
\end{algorithm}

Now we will analyze the seeds of $\tilde\Lambda$ depending on whether their order is an {\em old order}, i.e. at most $k-1$
or a {\em new order}, i.e., at least $k$.
The result that we will obtain at the end of this analysis is Theorem~\ref{t}. This theorem was already proved in \cite{seeds} from the perspective of the originary definition of seed introduced there.

\subsubsection*{Old-order seeds}

Suppose $p<k-1$ and let $\rho<\tilde u_{p}=u_p$.
We want to analyze whether $\lambda_s+1+\rho$ is an order-$p$ seed of $\tilde\Lambda$, that is, we want to determine whether $\tilde S_\ell=1$, for $\ell=\lambda_p+\rho$.
Since $p<k-1$,
we have $\ell+\lambda_s+1\leq \lambda_p+u_p+\lambda_s=\lambda_{p+1}+\lambda_s\leq \lambda_{k-1}+\lambda_s\leq c-2+\lambda_s\leq 2\lambda_s-2$. Hence, we are in one of the first three hypotheses of equation~\eqref{eq:seeds}.

{\em First case.} Suppose $0\leq \rho<u_p-(\lambda_s-c+1)$.
In this case
$\ell+\lambda_s+1< \lambda_{p+1}+c$.
On one hand this implies that $\ell+\lambda_s+1<2c$
and so, we are in one of the first two hypothesis of equation~\eqref{eq:seeds}.
On the other hand the same inequality implies that
$$\lambda_p+(\lambda_s+1)\leq \ell+\lambda_s+1<\lambda_{p+1}+c$$ and so
${\ell+\lambda_s+1}\not \in L+I_s$. Hence, $\tilde S_\ell=1$ if and only if $\Sigma_{\ell+\lambda_s+1}=1$, which is equivalent to $S_{\ell+\lambda_s-c+1}=1$.
  In this case we say we have an {\em old-order recycled seed} since this order-$p$ seed of $\tilde\Lambda$ is originated by an order-$p$ seed of $\Lambda$.
  Indeed, $\lambda_{p}\leq \ell+\lambda_s-c+1<\lambda_{p+1}$.
  Notice that if $\lambda_t=\lambda_s+1+\rho$ is an old-order recycled seed then $\lambda_t<\tilde c +u_p-1$.

{\em Second case.}
Suppose $u_p-(\lambda_s-c+1)\leq \rho\leq u_p-2$. In this case
$\ell+\lambda_s+1$ lies between
$\lambda_{p+1}+c$ and
$\lambda_{p+1}+\lambda_s-1$. Hence, 
$\ell+\lambda_s+1\in L+I_s$
and so $\tilde S_\ell=0$.

\smallskip

{\em Third case.}
Suppose $\rho=u_p-1$. In this case $\ell+\lambda_s+1=\lambda_{p+1}+\lambda_s$.

\vskip -.7cm\mbox{}
\begin{itemize}
\item
  If $u_{p+1}\leq \lambda_s-c$ then we need to distinguish whether $p=k-2$ or $p<k-2$.
  \begin{itemize}
    \item
  If $p=k-2$, then $\ell+\lambda_s+1=\lambda_{k-1}+\lambda_s\geq \lambda_{k-1}+u_{k-1}+c=2c$ while $\ell+\lambda_s+1=\lambda_{k-1}+\lambda_s\leq 2\lambda_s-2$, hence we are in the third hypothesis of equation~\eqref{eq:seeds} and, hence, $\tilde S_\ell=0$.
  \item If $p<k-2$, then 
$\ell+\lambda_s+1=\lambda_{p+2}+(\lambda_s-u_{p+1})$
    with $c\leq (\lambda_s-u_{p+1})\leq \lambda_s-1$, hence $\ell+\lambda_s+1\in L+I_s$ and we are in the first hypothesis of equation~\eqref{eq:seeds}. Then $\tilde S_\ell=0$.
    \end{itemize}
\item
  If $u_{p+1}>\lambda_s-c$ then $\lambda_{p+1}+\lambda_s<\lambda_{p+2}+c$, which implies that $\ell+\lambda_s+1<2c$ and so we are in the first or second hypothesis of equation~\eqref{eq:seeds}.
  If $p=k-2$ then $\ell+\lambda_s+1=\lambda_{k-1}+\lambda_s\not\in L+I_s$.
  If $p<k-2$ then by the previous inequality $\ell+\lambda_s+1=\lambda_{p+1}+\lambda_s<\lambda_{p+2}+c$ which implies, again, $\ell+\lambda_s+1\not\in L+I_s$.
  Hence, $\tilde S_\ell=1$ if and only if $\Sigma_{\ell+\lambda_s+1}=1$, which is equivalent to $S_{\ell+\lambda_s-c+1}=1$.
  In this case we say we have an {\em old-order new seed} since this order-$p$ seed of $\tilde\Lambda$ is originated by an order-$(p+1)$ seed of $\Lambda$.
  Indeed, $\lambda_{p+1}\leq \ell+\lambda_s-c+1<\lambda_{p+2}$.
\end{itemize}

Suppose now $p=k-1$ and let $\rho<\tilde u_{k-1}$. Notice that
$\tilde u_{k-1}=u_{k-1}$ if and only $\lambda_s\neq c$ and $\tilde u_{k-1}=u_{k-1}+1$ if $\lambda_s=c$.

{\em First case.} Suppose $0\leq \rho<u_{k-1}-(\lambda_s-c+1)$.
In this case
$\ell+\lambda_s+1< \lambda_{k}+c$.
On one hand this implies that $\ell+\lambda_s+1<2c$ and on the other hand,
since $\ell+\lambda_s+1\geq \lambda_{k-1}+\lambda_s+1$,
we have $\ell+\lambda_s+1\not\in L+I_s$.
Hence, $\tilde S_\ell=1$ if and only if $\Sigma_{\ell+\lambda_s+1}=1$, which is equivalent to $S_{\ell+\lambda_s-c+1}=1$.
  In this case we have again an {\em old-order recycled seed}.

{\em Second case.}
Suppose $u_{k-1}-(\lambda_s-c+1)\leq \rho< \tilde u_{k-1}$. In this case
$\ell+\lambda_s+1$ is at least
$\lambda_{k}+c=2c$.
Hence, $\tilde S_\ell=1$ if and only if $\ell+\lambda_s+1\geq 2\lambda_s-1$, that is, if and only if $\rho\geq \lambda_s-\lambda_{k-1}-2$, and so,
$$\lambda_s-\lambda_{k-1}-2\leq \rho\leq \tilde u_{k-1}-1.$$
This may occur only if
$\tilde u_{k-1}-1\geq \lambda_s-\lambda_{k-1}-2$, which is equivalent to
$\lambda_s\leq \tilde\lambda_k+1$ and which in turn is only possible if
\begin{itemize}
\item $\lambda_s=c+1$ and $\tilde\lambda_k=\lambda_k=c$, in which case $\tilde u_{k-1}=u_{k-1}$ and $u_{k-1}-1\leq \rho\leq u_{k-1}-1$.
\item $\lambda_s=c$ and $\tilde\lambda_k=\lambda_k+1$, in which case $\tilde u_{k-1}=u_{k-1}+1$ and $u_{k-1}-1\leq \rho\leq u_{k-1}$.
\end{itemize}
So, there are only old-order new-seeds $\lambda_s+1+\rho$ of order $k-1$ if and only if either $\lambda_s=c+1$ and $\rho=u_{k-1}-1$ or $\lambda_s=c$ and $\rho=u_{k-1}-1$ or $\rho=u_{k-1}$.

\subsubsection*{New-order seeds}

If $\lambda_s=c$ then all the possible seed orders of $\tilde \Lambda$ are already seed orders of $\Lambda$ and so, there are no new-order seeds.

Thus, we can assume that $\lambda_s\neq c$. Then $\tilde u_{s-1}=2$ and $\tilde u_{i}=1$ for any other new seed order $i$.

For $p=s-1$, since $\rho<\tilde u_p=\tilde u_{s-1}=2$ we only need to consider the cases $\rho=0$ and $\rho=1$. But either for $\rho=0$ and for $\rho=1$ the non-gaps $\lambda_s+1+\rho$ correspond to order-$(s-1)$ seeds. Indeed the corresponding $\ell$ is $\ell=\lambda_{s-1}+\rho=\lambda_s+\rho-1$ and hence, $\ell+\lambda_s+1=2\lambda_s+\rho$, which in both cases lies between
$2\lambda_s-1$ and $2\lambda_s+1$. There are no more order-$(s-1)$ seeds since $\tilde u_{s-1}=2$.

For $p=s-2$, if $p$ is indeed a new seed order, then $\tilde u_{s-2}=1$. 
Then $\rho$ can only be $0$ 
Now, $\ell+\lambda_s+1=\lambda_{s-2}+\lambda_s+1=2\lambda_s-1$ and by equation~\eqref{eq:seeds}, $\tilde S_\ell=1$ and, so, $\lambda_s+1$ is the unique order-$(s-2)$ seed.

For $k\leq p< s-2$, necessarily $\tilde u_{p}=1$ and, so, the unique possible seed would be $\rho=0$. In this case $\ell=\lambda_p$ and $\ell+\lambda_s+1=\lambda_p+\lambda_s+1<2\lambda_s-1$. Hence,
by equation~\eqref{eq:seeds}, $\tilde S_\ell=0$ and, so, $\lambda_s+1+\rho$ is not an order-$(s-2)$ seed.

The previous results can be summarized in the next theorem.

\begin{theorem}\label{t}
  Suppose that $\Lambda$ is a semigroup of rank $k$ and conductor $c$ and let $\tilde \Lambda=\Lambda\setminus\{\lambda_s\}$ with $s\geq k$.
The seeds of $\tilde \Lambda$ of order $p$ are exactly the ones in the next list.

\begin{enumerate}
\item If $p<k$ ({\bf Old-order seeds}):
    \begin{itemize}
  \item {\bf Old-order recycled seeds:}
    \begin{itemize} \item Any order-$p$ seed $\lambda_t$ of $\Lambda$ with \,$t>s$ (necessarily with $\lambda_t-c(\Lambda\setminus\{\lambda_s\})<u_p-1$) \end{itemize}
\item {\bf Old-order new seeds:}
  \begin{itemize}
  \item $\lambda_t=\lambda_s+u_p$\quad\quad\quad\quad\quad\quad\quad\quad\quad\quad{if $p<k-1$ and $\lambda_s$ is an order-$(p+1)$ seed of $\Lambda$}
  \item $\lambda_t=\lambda_s+u_p$ and $\lambda_t=\lambda_s+u_p+1$\quad{if $p=k-1$ and $\lambda_s=c$}
  \item $\lambda_t=\lambda_s+u_p$\quad\quad\quad\quad\quad\quad\quad\quad\quad\quad{if $p=k-1$ and $\lambda_s=c+1$}
\end{itemize}
\end{itemize}
\item If $p\geq k$ ({\bf New-order seeds})
\begin{itemize}
\item $\lambda_s+1$ if $p=s-2$
\item $\lambda_s+1$ and $\lambda_s+2$ if $p=s-1$
\end{itemize}
\end{enumerate}

\end{theorem}

\section{New generators, new strong generators and newly strong generators}
\label{s:specialcases}

Let $\Lambda$ be a numerical semigroup with $c(\Lambda)\geq \lambda_1$. A {\em new generator} of $\Lambda$ is a generator of $\Lambda$ that is not a generator of $\Lambda\cup F(\Lambda)$.
A {\em new strong generator} of $\Lambda$ is a strong generator of $\Lambda$ that is not a generator of $\Lambda\cup F(\Lambda)$.
A {\em newly strong generator} of $\Lambda$ is a generator of $\Lambda$ and also a generator of $\Lambda\cup F(\Lambda)$, that is weak in $\Lambda\cup F(\Lambda)$ and strong in $\Lambda$. 

Recall that order-$i$ seeds are elements between $c(\Lambda)$ and $c(\Lambda)+u_i(\Lambda)-1$ and they are represented by $1$'s in the corresponding positions of $S(\Lambda)$ between positions $\lambda_i$ and $\lambda_i+u_i(\Lambda)-1=\lambda_{i+1}-1$. Recall also that right generators of $\Lambda$ are exactly the order-$0$ seeds of $\Lambda$ and that a right generator is strong if and only if it is an order-$1$ seed. In particular, right generators lie between $c(\Lambda)$ and $c(\Lambda)+m(\Lambda)-1$ while strong generators lie between $c(\Lambda)$ and $c(\Lambda)+u(\Lambda)-1$, where $u(\Lambda)=u_1(\Lambda)$.

For ordinary semigroups (i.e. semigroups of rank $k=1$), the notion of strong generator is not defined. But we can define the {\em ordinarily strong generators} of an ordinary semigroup as those generators than when are taken away, the obtained semigroup has new generators. One can check that the ordinary semigroup $\Lambda_m$ of multiplicity $m$ has exactly two ordinarily strong generators: $m$ and $m+1$. The semigroup $\Lambda_m\setminus\{m\}$ has two new generators, $2m$ and $2m+1$; the semigroup $\Lambda_m\setminus\{m+1\}$ has one new generator, $2m+1$.

Let $\Lambda$ be a numerical semigroup with $S=S(\Lambda)$, $c=c(\Lambda)$, $m=m(\Lambda)$, $u=u(\Lambda)=u_(\Lambda)$, and $v=v(\Lambda)=u_2(\Lambda)$.

\begin{enumerate}[{\bf (1)}]
\item Existence of right generators. The element $\mu\geq c$ is a right generator of $\Lambda$ if and only if $S_{\mu-c}=1$. The number of right generators of $\Lambda$ is $w_0^{m-1}(S)$. By Lemma~\ref{l:lowrank}, in the particular case $k(\Lambda)=1$, we have $w_0^{m-1}(S)=m$ and in the particular case $k(\Lambda)=2$, we have $w_0^{m-1}(S)=m-1$
\item Existence of strong generators. If $k(\Lambda)=1$ there are no strong generators by definition but there are two ordinarily strong generators. If $k(\Lambda)>1$, the element $\mu\geq c$ is a strong generator of $\Lambda$ if and only if $S_{\mu-c}=S_{\mu+m-c}=1$. The number of strong generators of $\Lambda$ is $w_0^{u-1}(S\land (S\ll m))$.
  In the particular case $k(\Lambda)=2$, by Lemma~\ref{l:lowrank}, $w_0^{u-1}(S\land (S\ll m))=\left\{\begin{array}{ll}u-1&\mbox{ if }m-u<u\\u&\mbox{ otherwise.}\end{array}\right.$
\item Existence of new generators. Suppose $\mu_1$ is a right generator of $\Lambda$.
  \begin{itemize}
    \item 
      Suppose $k(\Lambda)=1$ then $\Lambda$ is ordinary and, as discussed above, $\Lambda\setminus\{\mu_1\}$ has exactly two new generators if $\mu_1=m$, one new generator if $\mu_1=m+1$ and zero new generators otherwise.
    \item 
      If $k(\Lambda)>1$, by Theorem~\ref{t} (old-order new seeds, $p<k-1$) there is at most one new generator of $\Lambda\setminus\{\mu_1\}$, which must be $\mu_1+m$. Indeed, $\mu_1+m$ is a new generator if and only $\mu_1$ is a strong generator of $\Lambda$ and so, if and only if $S_{\mu_1-c}=S_{\mu_1+m-c}=1$. This situation occurs $w_0^{u-1}(S\land(S\ll m))$ times.
  \end{itemize}

\item Existence of new strong generators. Suppose $\mu_1$ is a right generator of $\Lambda$.
  \begin{itemize}
  \item
    Suppose $k(\Lambda)=1$.
    \begin{itemize}
      \item If $\mu_1=c$ then $k(\Lambda\setminus\{\mu_1\})=1$
    and $\Lambda\setminus\{\mu_1\}$ has no strong generators since in this case there are no order-$1$ seeds.
  \item 
    If $\mu_1=c+1$, we just saw that there is one new generator of $\Lambda\setminus\{\mu_1\}$, which is $\mu_1+m$.
    By Theorem~\ref{t} (new-order new seeds, $p=s-1$)
    the new generator is strong if and only if $\mu_1+m=\mu_1+1$ or $\mu_1+m=\mu_1+2$. In the first case we have a contradiction with $k(\Lambda)=1$. So, there exists one new strong generator if $m=2$ for $\mu_1=3$, and zero new strong generators otherwise.

  \item If $\mu_1>c+1$, we already saw that there are no new generators of $\Lambda\setminus\{\mu_1\}$.
    \end{itemize}
  \item 
    If $k(\Lambda)=2$
    we already saw that
    there is one new generator of $\Lambda\setminus\{\mu_1\}$, which is $\mu_2=\mu_1+m$, only if $\mu_1$ is a strong generator of $\Lambda$.
By Lemma~\ref{l:lowrank} for the case of rank $k=2$, this is equivalent to $\mu_1\leq c+u-1$. 
    Now, $\mu_2$ may be strong in $\Lambda\setminus\{\mu_1\}$
    if and only if it is a seed of order $1$ of $\Lambda\setminus\{\mu_1\}$. We know that order-$1$ seeds of $\Lambda\setminus\{\mu_1\}$ are at most
    $\mu_1+u_1(\Lambda\setminus\{\mu_1\})$. In particular, we need $\mu_2=\mu_1+m\leq \mu_1+u_1(\Lambda\setminus\{\mu_1\})$ and so $u_1(\Lambda\setminus\{\mu_1\})\geq m$ which implies $u_1(\Lambda\setminus\{\mu_1\})=m$. This in turn implies that either $u=m-1$ and $\mu_1=c$ or $u=m$ and $\mu_1>c$.
    Now, since $\mu_2=\mu_1+m$, it can not be an (old-)order-$1$ recycled seed of $\Lambda\setminus\{\mu_1\}$. Then it needs to be an (old-)order-$1$ new seed of $\Lambda\setminus\{\mu_1\}$.
    Then, by Theorem~\ref{t} (for old-order new seeds with $p=k-1$), $\mu_2$ is a strong generator of $\Lambda\setminus\{\mu_1\}$ if and only if either
    \begin{itemize}
      \item
        $u=m-1$ and $\mu_1=c$
      \item $u=m$ and $\mu_1=c+1$.
    \end{itemize}
    So, there exists one new strong generator if $u=m-1$ for $\mu_1=c$ and there exists one new strong generator if $u=m$ for $\mu_1=c+1$. There exist zero new strong generators otherwise.
  \item
  Suppose $k(\Lambda)\geq 3$.
  A new generator $\mu_2$ of $\Lambda\setminus\{\mu_1\}$ may only be strong if $u=m$. Indeed, first of all,
  $u_0(\Lambda)=u_0(\Lambda\setminus\{\mu_1\})$ and
  $u_1(\Lambda)=u_1(\Lambda\setminus\{\mu_1\})$.
  Second, in order for $\mu_2$ to be a new generator, $\mu_2=\mu_1+m$. Third, for $\mu_2$ to be strong in $\Lambda\setminus\{\mu_1\}$ it must be an order-$1$ seed of $\Lambda\setminus\{\mu_1\}$. Since the conductor of $\Lambda\setminus\{\mu_1\}$ is $\mu_1+1$, an order-$1$ seed of $\Lambda\setminus\{\mu_1\}$ is at most $\mu_1+u$. Hence $\mu_1+m\leq \mu_1+u$ and this is only possible if $u=m$.
      Now, since $\mu_2=\mu_1+m$, it can not be an (old-)order-$1$ recycled seed of $\Lambda\setminus\{\mu_1\}$. Then it needs to be an (old-)order-$1$ new seed of $\Lambda\setminus\{\mu_1\}$.
      In this case, by Theorem~\ref{t} (old-order new seeds with $p<k-1$), $\mu_1$ needs to be an order-$2$ seed of $\Lambda$ and so, $S_{\mu_1-c+\lambda_2}=S_{\mu_1-c+2m}=1$. Hence, there exists one new strong generator of $\Lambda\setminus\{\mu_1\}$ if $u=m$ and $S_{\mu_1-c}=S_{\mu_1-c+m}=S_{\mu_1-c+2m}=1$.
This situations occurs $w_0^{v-1}(S\land(S\ll u)\land(S\ll(m+u)))$ times whenever $u=m$. Otherwise there exist no new strong generators.
  \end{itemize}

  \item Existence of newly strong generators. If $\mu_2$ is a newly strong generator of $\Lambda\setminus\{\mu_1\}$, then, on one hand $c\leq \mu_1<\mu_2<c+m$ and so $\mu_2-\mu_1<m$ and, on the other hand $\mu_2$ is an (old-)order-$1$ new seed of $\Lambda\setminus\{\mu_1\}$. 
    \begin{itemize}
    \item If $k(\Lambda)=1$, by Lemma~\ref{l:lowrank} $\mu_1,\mu_2$ with $\mu_1<\mu_2$ are right generators if and only if $m\leq \mu_1<\mu_2<2m$.
      By Theorem~\ref{t} (new-order new seeds with $p=k$), $\mu_2$ is a newly strong generator of $\Lambda\setminus\{\mu_1\}$ if and only if one of the following options holds.
      \begin{itemize}
      \item
        $\mu_1=\lambda_3=m+2$ and $\mu_2=\mu_1+1=m+3$ and $m\geq 4$ 
      \item
        $\mu_1=\lambda_2=m+1$ and $\mu_2=\mu_1+1=m+2$ and $m\geq 3$
      \item
        $\mu_1=\lambda_2=m+1$ and $\mu_2=\mu_1+2=m+3$ and $m\geq 4$
      \end{itemize}
      So, if $m\leq 2$ there exist zero newly strong generators for any $\mu_1$;
      if $m=3$ there exists one newly strong generator for $\mu_1=m+1$ and zero newly strong generators otherwise; if $m\geq 4$ there exist two newly strong generators for $\mu_1=m+1$, one newly strong generator for $\mu_1=m+2$,
      and zero newly strong generators otherwise.

      We add that there exists one newly ordinarily strong generator for $\mu_1=m$, which is $m+2$.
    \item If $k(\Lambda)=2$, by Theorem~\ref{t} (old-order new seeds with $p=k-1$), one of the following options must hold.
      \begin{itemize}
      \item 
        $\mu_1=c=m+u$ and $\mu_2=\mu_1+u=m+2u$ with $\mu_1$, $\mu_2$ right generators of $\Lambda$. That is, by Lemma~\ref{l:lowrank}, $\mu_1-c\neq m-u$ and $\mu_2-c\neq m-u$, i.e., $u\neq m$ and $m\neq 2u$.
        Conversely, one can check that if $u<m$ and $m\neq 2u$ then $m+2u$ is a newly strong generator of $\Lambda\setminus\{m+u\}=\{0,m,m+u+1,m+u+2,\dots\}$.
      \item
        $\mu_1=c=m+u$ and $\mu_2=\mu_1+u+1=m+2u+1$ with $\mu_1$, $\mu_2$ right generators of $\Lambda$. That is, by Lemma~\ref{l:lowrank}, $\mu_1-c\neq m-u$ and $\mu_2-c\neq m-u$, i.e., $u\neq m$ and $m\neq 2u+1$, which together with the condition $\mu_2<c+m$ results in the conditions $u<m-1$ and $m\neq 2u+1$.
        Conversely, one can check that if $u<m-1$ and $m\neq 2u+1$ then $m+2u+1$ is a newly strong generator of $\Lambda\setminus\{m+u\}=\{0,m,m+u+1,m+u+2,\dots\}$.
      \item
        $\mu_1=c+1=m+u+1$ and $\mu_2=\mu_1+u=m+2u+1$ with $\mu_1$, $\mu_2$ right generators of $\Lambda$. That is, by Lemma~\ref{l:lowrank}, $\mu_1-c\neq m-u$ and $\mu_2-c\neq m-u$, i.e., $u\neq m-1$ and $m\neq 2u+1$, which together with the condition $\mu_2-\mu_1<m$ results in the conditions $u< m-1$ and $m\neq 2u+1$.
        Conversely, one can check that if $u<m-1$ and $m\neq 2u+1$ then $m+2u+1$ is a newly strong generator of $\Lambda\setminus\{m+u+1\}=\{0,m,m+u,m+u+2,m+u+3,\dots\}$.
      \end{itemize}
      So, if $u<m$ and $m\neq 2u$ there exists one newly strong generator for $\mu_1=m+u$, while if $u<m-1$ and $m\neq 2u+1$ there exists one newly strong generator for $\mu_1=m+u$ (different than the one just mentioned) and a newly strong generator for $\mu_1=m+u+1$. There exist zero newly strong generators otherwise.
    \item If $k(\Lambda)\geq 3$, by Theorem~\ref{t} (old-order new seeds with $p<k-1$), $\mu_2=\mu_1+u$ and $\mu_1$ needs to be an order-$2$ seed of $\Lambda$ and so, $S_{\mu_1-c+\lambda_2}=S_{\mu_1-c+m+u}=1$. Furthermore, since $\mu_2-\mu_1<m$ and $\mu_2=\mu_1+u$, it follows that $u<m$.
    Hence, there exists one newly strong generator of $\Lambda\setminus\{\mu_1\}$ if $u<m$ and $S_{\mu_1-c}=S_{\mu_1-c+u}=S_{\mu_1-c+u+m}=1$. 
This situation occurs $w_0^{v-1}(S\land(S\ll u)\land(S\ll(m+u)))$ times whenever $u<m$. Otherwise there exist no newly strong generators. 
    \end{itemize}

\end{enumerate}

\section{The great-grandchildren of a numerical semigroup}
\label{s:greatgrandchildren}

We explain now how to obtain the set of all great-grandchildren of a semigroup in the semigroup tree just from its seeds, and the parameters $m$, $u$, $v$.


Each particular case of the lemma is illustrated with an example at the appendix. In those examples we represent the tree of descendants of semigroups with different characteristics.
Each semigroup in the tree is represented by its table of seeds.

\begin{theorem}\label{l:ggc}Let $S=S(\Lambda)$, $m=m(\Lambda)$, $u=u(\Lambda)$, $v=v(\Lambda)$.
  
  The number $n_c$ of children, the number $n_{gc}$ of grandchildren, and the number $n_{ggc}$ of great-grandchildren of $\Lambda$ is, respectively,

\begin{eqnarray*}
n_c(\Lambda)&=&\left\{\begin{array}{lll}
1 &\mbox{ if } k(\Lambda)= 0, &\text{ (see Fig.~\ref{fig:exa})}\\
m &\mbox{ if } k(\Lambda)= 1, &\text{ (see Fig.~\ref{fig:exa})}\\
m-1 &\mbox{ if } k(\Lambda)= 2, &\text{ (see Fig.~\ref{fig:exa})}\\
w_0^{m-1}(S)&\mbox{ if } k(\Lambda)\geq 3, &\text{ (see Fig.~\ref{fig:exb}-Fig.~\ref{fig:exf})}\\
\end{array}\right.
\\
n_{gc}(\Lambda)&=&\left\{\begin{array}{lll}
2 &\mbox{ if } k(\Lambda)= 0, &\text{(see Fig.~\ref{fig:exa})}\\
\binom{m}{2} + 3&\mbox{ if }k(\Lambda)= 1, &\text{ (see Fig.~\ref{fig:exb}-Fig.~\ref{fig:exd})}\\
\binom{n_c(\Lambda)}{2}+w_0^{u-1}(S\land(S\ll m))&\mbox{ if }k(\Lambda)>1, &\text{ (see Fig.~\ref{fig:exe}-Fig.~\ref{fig:exf})}\\\end{array}\right.
\\
n_{ggc}(\Lambda)&=&
  \left\{\begin{array}{lll}
  4 &\mbox{ if } k(\Lambda)= 0, &\text{ (see Fig.~\ref{fig:exa})}\\
\binom{m}{3}+3m+1
&\mbox{ if } k(\Lambda)=1, m=2,3, &\text{ (see Fig.~\ref{fig:exb} and Fig~\ref{fig:exc})}\\
\binom{m}{3}+3m+3
&\mbox{ if } k(\Lambda)=1, m\geq 4, &\text{ (see Fig.~\ref{fig:exd})}\\
  \binom{m-1}{3}
+(u-\delta_a)(m-2)+\delta_b+2\delta_{c}+\delta_d
&\mbox{ if } k(\Lambda)=2, &\text{ (see Fig.~\ref{fig:exe})}\\
\binom{n_c(\Lambda)}{3}
+w_0^{u-1}(S\land(S\ll m))(n_c(\Lambda)-1)\\
\phantom{\binom{n_c(\Lambda)}{3}}+w_0^{v-1}(S\land(S\ll u)\land(S\ll(u+m)))
&\mbox{ if } k(\Lambda)\geq3, &\text{ (see Fig.~\ref{fig:exf})}\\
  \end{array}\right.
\\
  \end{eqnarray*}
  where
\begin{eqnarray*}
  \delta_a&=&\left\{\begin{array}{ll}1&\mbox{ if }m<2u\\0&\mbox{ otherwise}\end{array}\right.\\
  \delta_{b}&=&\left\{\begin{array}{ll}1&\mbox{ if }u=m\mbox{ or }2u\neq m  \\0&\mbox{ otherwise}\end{array}\right.\\
  \delta_c&=&\left\{\begin{array}{ll}1&\mbox{ if }u<m-1\mbox{ and }2u+1\neq m  \\0&\mbox{ otherwise}\end{array}\right.\\
  \delta_d&=&\left\{\begin{array}{ll}1&\mbox{ if }u=m-1\\0&\mbox{ otherwise}\end{array}\right.\\
\end{eqnarray*}

\end{theorem}

\begin{proof}
  The unique child of ${\mathbb N}$ is ${\mathbb N}\setminus \{1\}$,
  while for a semigroup $\Lambda\neq {\mathbb N}$ the number of children is the number of its right generators. Then the formula for $n_{c}$ follows from {\bf (1)}.

  Let us prove now the formula for $n_{gc}(\Lambda)$.
  The number of grandchildren of ${\mathbb N}$ is $2$ since this is the number of semigroups of genus $2$.

  If $\Lambda\neq {\mathbb N}$, then a grandchild of $\Lambda$ is a semigroup of the form $\Lambda\setminus\{\mu_1,\mu_2\}$ with $c\leq \mu_1<\mu_2$ and either
\begin{enumerate}[(a)]
  \item $\mu_1,\mu_2$ are right generators of $\Lambda$,
  \item $\mu_1$ is a right generators of $\Lambda$, and $\mu_2$ is a new generator of $\Lambda\setminus\{\mu_1\}$.
  \end{enumerate}
  
  Case (a) can be obtained in $\binom{n_c(\Lambda)}{2}$ different ways.
  Case (b), for $k(\Lambda)=1$ can be obtained in three different ways, taking $\mu_1$ as one of the two ordinarily strong generators and $\mu_2$ as one of the new generators corresponding to $\mu_1$ ({\bf (3)}).
Case (b), for $k(\Lambda)\geq 2$ can be done as many times as the number of strong generators of $\Lambda$, which is 
$w_0^{u-1}(S\land(S\ll m))$ ({\bf (2)}).

Let us prove now the formula for $n_{ggc}$.
If $\Lambda={\mathbb N}$ then we know that the number of its great-grandchildren is $4$. It can be just observed in Figure~\ref{f:tree}.

If $\Lambda\neq {\mathbb N}$, then a great-grandchild of $\Lambda$ is a semigroup of the form $\Lambda\setminus\{\mu_1,\mu_2,\mu_3\}$ with $c\leq \mu_1<\mu_2<\mu_3$ and either
  \begin{enumerate}
  \item[(a)] $\mu_1,\mu_2,\mu_3$ are right generators of $\Lambda$,
  \item[(b)] $\mu_1,\mu_2$ are right generators of $\Lambda$, and $\mu_3$ is a new generator of $\Lambda\setminus\{\mu_1\}$ or a new generator of $\Lambda\setminus\{\mu_2\}$,
  \item[(c)] $\mu_1$ is a right generator of $\Lambda$ and $\mu_2,\mu_3$ are new generators of $\Lambda\setminus\{\mu_1\}$,
  \item[(d)] $\mu_1,\mu_2$ are right generators of $\Lambda$, $\mu_3$ is not a right generator of $\Lambda$, nor a right generator of $\Lambda\setminus\{\mu_1\}$, nor a right generator of $\Lambda\setminus\{\mu_2\}$, but it is a new generator of $\Lambda\setminus\{\mu_1,\mu_2\}$,
  \item[(e)] 
$\mu_1$ is a right generator of $\Lambda$, $\mu_2$ is a new generator of $\Lambda\setminus\{\mu_1\}$ and $\mu_3$ is a new generator of $\Lambda\setminus\{\mu_1,\mu_2\}$.
  \end{enumerate}

  Now we are going to analyze each case separately.
  \begin{itemize}
    \item   The situation (a) is counted by
  $\binom{n_c(\Lambda)}{3}$.
      In particular, for
      $k(\Lambda)=1$ this equals $\binom{m}{3}$ and for $k(\Lambda)=2$ this equals $\binom{m-1}{3}$.
        \begin{center}\begin{tabular}{|c|c|c|}
            \hline
            $k(\Lambda)=1$& $k(\Lambda)=2$ & $k(\Lambda)\geq 3$\\\hline
$\binom{m}{3}$ & $\binom{m-1}{3}$ & $\binom{n_c(\Lambda)}{3}$\\\hline
        \end{tabular}\end{center}
      
    \item
      Suppose now situation (b). For $k(\Lambda)=1$, considering the two ordinarily strong generatrors, this may occur in the following cases:
      \begin{itemize}
      \item $\mu_1=m$, $\mu_2=m+1$, $\mu_3=2m$
      \item $\mu_1=m$, $\mu_2=m+1$, $\mu_3=2m+1$
      \item $\mu_1=m$, $\mu_2\in\{m+2,\dots,2m-1\}$, $\mu_3=2m$
      \item $\mu_1=m$, $\mu_2\in\{m+2,\dots,2m-1\}$, $\mu_3=2m+1$
      \item $\mu_1=m+1$, $\mu_2\in\{m+2,\dots,2m-1\}$, $\mu_3=2m+1$
      \end{itemize}
      So, this situation occurs $2+3(m-2)=3m-4$ times.
      For $k(\Lambda)\geq 2$, the situation in (b) corresponds to the number of strong generators multiplied by the number of right generators minus one.
      By {\bf (1)} and {\bf (2)}, this equals $w_0^{u-1}(S\land(S\ll m))(n_c(\Lambda)-1)$, which for the particular case $k(\Lambda)=2$ is $(u-\delta_a)(m-2)$.

              \begin{center}\begin{tabular}{|c|c|c|}
            \hline
            $k(\Lambda)=1$& $k(\Lambda)=2$ & $k(\Lambda)\geq 3$\\\hline
$3m-4$ & $(u-\delta_a)(m-2)$ &  $w_0^{u-1}(S\land(S\ll m))(n_c(\Lambda)-1)$\\\hline
        \end{tabular}\end{center}

            \item Situation (c) only occurs once for $k(\Lambda)=1$, for $\mu_1=m,\mu_2=2m,\mu_3=2m+1$.

              \begin{center}\begin{tabular}{|c|c|c|}
            \hline
            $k(\Lambda)=1$& $k(\Lambda)=2$ & $k(\Lambda)\geq 3$\\\hline
            $1$  & & \\\hline
        \end{tabular}\end{center}
      \item   
        Situation (d) for $k(\Lambda)>1$ or for $k(\Lambda)=1$ with $\mu_1>m$, is equivalent to $\mu_1,\mu_2$ being right generators of $\Lambda$, $\mu_2$ being a newly strong generator of $\Lambda\setminus\{\mu_1\}$ and $\mu_3$ being the unique new generator of $\Lambda\setminus\{\mu_1,\mu_2\}$.
        According to {\bf (5)} this occurs as many times as indicated in the table
        \begin{center}\begin{tabular}{|c|c|c|}
            \hline
            $k(\Lambda)=1$& $k(\Lambda)=2$ & $k(\Lambda)\geq 3$\\\hline
$0$ if $m\leq 2$ & $\delta_{b'}+2\delta_c$ & 
            $(1-\delta_{b''})w_0^{v-1}(S\land(S\ll u)\land(S\ll(u+m)))$\\
            $1$ if $m=3$ & & \\
            $3$ if $m\geq 4$ & & \\\hline
        \end{tabular}\end{center}
where
      \begin{eqnarray*}
    \delta_{b'}&=&\left\{\begin{array}{ll}1&\mbox{ if }u<m\mbox{ and }2u\neq m  \\0&\mbox{ otherwise}\end{array}\right.
  \end{eqnarray*}
    \begin{eqnarray*}
  \delta_{b''}&=&\left\{\begin{array}{ll}1&\mbox{ if }u=m  \\0&\mbox{ otherwise}\end{array}\right.\\
\end{eqnarray*}

    Situation (d) for $k(\Lambda)=1$ and $\mu_1=m$ may only occur for $\mu_2$ equal to one of the ordinarily strong generators of $\Lambda\setminus\{m\}$, which are $m+1$ and $m+2$.
    For $m=2$,
    the cases in which this occurs are excatly when either $(\mu_1,\mu_2,\mu_3)=(2,3,6)$ or $(\mu_1,\mu_2,\mu_3)=(2,3,7)$. For $m>2$, the cases in which this occurs are excatly when either $(\mu_1,\mu_2,\mu_3)=(m,m+1,2m+2)$, $(\mu_1,\mu_2,\mu_3)=(m,m+1,2m+3)$, or $(\mu_1,\mu_2,\mu_3)=(m,m+2,2m+2)$.
So, we need to add to the previous table these cases.
\begin{center}\begin{tabular}{|c|c|c|}
            \hline
            $k(\Lambda)=1$& $k(\Lambda)=2$ & $k(\Lambda)\geq 3$\\\hline
            $2$ if $m=2$ & & \\
            $3$ if $m\geq 3$ & & \\\hline
        \end{tabular}\end{center}

\item
Situation (e) is equivalent to $\mu_1$ being a right generator of $\Lambda$, $\mu_2$ being a new strong generator of $\Lambda\setminus\{\mu_1\}$ and $\mu_3$ being a new generator of $\Lambda\setminus\{\mu_1,\mu_2\}$.
    For $m=2$ this occurs exactly for $(\mu_1,\mu_2,\mu_3)=(2,4,7)$ and for $(\mu_1,\mu_2,\mu_3)=(3,5,7)$.
Otherwise, according to {\bf (4)} this occurs as many times as indicated in the table.

        \begin{center}\begin{tabular}{|c|c|c|}
            \hline
            $k(\Lambda)=1$& $k(\Lambda)=2$ & $k(\Lambda)\geq 3$\\\hline
$2$ if $m= 2$ & $\delta_{b''}+\delta_d$ & 
            $\delta_{b''}\cdot w_0^{v-1}(S\land(S\ll u)\land(S\ll(u+m)))$ \\
            $0$ if $m\geq 3$ & & \\\hline
        \end{tabular}\end{center}
  \end{itemize}

The result follows from the fact that $\delta_b=\delta_{b'}+\delta_{b''}$.
  
\end{proof}

\begin{example}
  For instance, the semigroup $\Lambda=\{0,3,6,8,9,10,\dots\}$ (with $S(\Lambda)=10110111$) satisfies
  \begin{itemize}
  \item $n_c(\Lambda)=2$, since $w_0^{m-1}(S)=w_0^2(S)=w(101)=2$.
  \item $n_{gc}(\Lambda)=3$, since $\binom{n_c(\Lambda)}{2}+w_0^{u-1}(S\land(S\ll m))=1+w_0^{3}(101\land 101)=1+2=3$.
  \item $n_{ggc}(\Lambda)=3$, since $m=u=3$, $v=2$,
    $w_0^2(S\land(S\ll m))=w(101\land 101)=2$, 
and $w_0^{1}(S\land(S\ll u)\land (S\ll(u+m))=w(10\land 10\land 11)=1$, 
the formula gives $\binom{2}{3}+2(2-1)+1=3$.
  \end{itemize}
\end{example}

\section{The seeds algorithm revisited}
\label{s:performance}


The previous results can be used to define a revisited version of the seeds algorithm whose implementation in C++ can be seen in the appendix. The bit count required by Lemma~\ref{l:ggc} is performed using Brian Kernighan's algorithm.
Parallelizing the implementation gives very good computation times.
A parallel implementation of the algorithm can be found in
\\\phantom{mmmmmmmmmm}{\tt https://github.com/mbrasamoros/seeds-algorithm}.\\
The performance of the algorithm is next compared to that of the Fromentin-Hivert algorithm \cite{FromentinHivert,FH-code} and to that of the RGD algorithm \cite{rgd,RGD-code}.
Using 8 workers in all algorithms we obtain the following table, where computation time is measured in seconds.

\begin{center}
\begin{tabular}{|l| r r r r r r r r r r | }
\hline
genus  & 40 & 42 & 44 & 46 & 48 & 50 & 52 & 54 & 56 & 58 \\
\hline
FH  & 1 & 3 & 10 & 28 & 75 & 207 & 537 & 1410 & 3977 & 10620 \\
RGD  & 2 & 4 & 9 & 23 & 62 & 158 & 403 & 1064 & 2892 & 7462 \\
Seeds & 1 & 2 & 6 & 15 & 41 & 107 & 284 & 758 & 1962 & 5257 \\
\hline
\end{tabular}
\end{center}
Using 12 workers in all algorithms we obtain the next table.
\begin{center}
\begin{tabular}{|l| r r r r r r r r r r r| }
\hline
genus  & 40 & 42 & 44 & 46 & 48 & 50 & 52 & 54 & 56 & 58 & 60 \\
\hline
FH  & 1 & 2 & 7 & 19 & 53 & 145 & 372 & 978 & 2760 & 7398 & 21880 \\
  RGD  & 1 & 3 & 6 & 18 & 45 & 121 & 291 & 799 & 2101 & 5292 & 13785 \\
Seeds  & 1 & 2 & 4 & 11 & 27 & 73 & 195 & 503 & 1329 & 3556 & 9459 \\
\hline
\end{tabular}
\end{center}
And using 24 workers we obtain this other table.
\begin{center}
\begin{tabular}{|l| r r r r r r r r r r r| }
\hline
genus  & 40 & 42 & 44 & 46 & 48 & 50 & 52 & 54 & 56 & 58 & 60 \\
\hline
FH  & 1 & 2 & 6 & 17 & 45 & 126 & 332 & 876 & 2494 & 6600 & 19090 \\
  RGD  & 1 & 2 & 5 & 13 & 38 & 88 & 244 & 625 & 1551 & 4458 & 10382 \\
Seeds  & 1 & 1 & 4 & 9 & 25 & 67 & 179 & 462 & 1226 & 3307 & 8586 \\
\hline
\end{tabular}
\end{center}

\mut{
\section{GUMs i GUMVs?}

In the implementation of the RGD algorithm, as explained in \cite{rgd}, the exploration of the semigroup tree is divided among the subtrees of the semigroups of a given multiplicity and jump.

In the implementation of the seeds algorithm here presented,
the exploration of the semigroup tree is divided among the subtrees of the semigroups of a given multiplicity, jump, and second jump.

Hence, for a fixed genus $g$, it will be interesting to have an estimate of the behavior of the parameters
\begin{eqnarray*}
  n_{g,m,u}&=&\#\{\Lambda: \Lambda \mbox{ is a numerical semigroup of genus } g, \mbox{ with }m(\Lambda)=m,u(\Lambda)=u\}\\
  n_{g,m,u,v}&=&\#\{\Lambda: \Lambda \mbox{ is a numerical semigroup of genus } g, \mbox{ with }m(\Lambda)=m,u(\Lambda)=u,v(\Lambda)=v\}\\ 
  \end{eqnarray*}for different values of $m,u,v$.
}

\section{Computation of Eliahou semigroups}
\label{s:eliahou}

Fix a numerical semigroup $\Lambda$.
Let $k=\rank(\Lambda)$, let $p$ be the number of primitive elements of $\Lambda$, and let $r$ the number of right generators of $\Lambda$.
The Wilf conjecture states that $c(\Lambda)\leq kp$ \cite{Wilf}.
Let $q=\left\lceil\frac{c(\Lambda)}{m(\Lambda)}\right\rceil$ and $\rho=qm(\Lambda)-c(\Lambda)$.
We say that $\Lambda$ is an Eliahou semigroup if the constant $E(\Lambda)=k(p-r)-q(m-r)+\rho$ is negative.
Shalom Eliahou proved that if the Wilf conjecture does not hold for a semigroup, then it must be an Eliahou semigroup \cite{eliahou}.
Hence, verifying that the Eliahou semigroups up to a certain genus satisfy the Wilf conjecture proves the Wilf conjecture for all semigroups up to this genus.
Eliahou semigroups are very rare. Let $\langle a,b,c\rangle\mid_{\kappa}$ be the minimum semigroup containing $a,b,c$ and all integers larger than or equal to $\kappa$.
Jean Fromentin found that the unique Eliahou semigroups of genus $g\leq 60$ are exactly
\begin{eqnarray*}
    \varepsilon_1&=&\langle 14,22,23\rangle\mid_{56},\\
   \varepsilon_2&=&\langle 16,25,26\rangle\mid_{64},\\
 \varepsilon_3&=&\langle 17,26,28\rangle\mid_{68},\\
    \varepsilon_4&=&\langle 17,27,28\rangle\mid_{68},\\
    \varepsilon_5&=&\langle 18,28,29\rangle\mid_{72},\\
\end{eqnarray*}
 which have genus $43, 51, 55, 55$ and $59$, respectively.
In \cite{mdai} it was computed that the unique Eliahou semigroups with genus between $61$ and $65$ are exactly
\begin{eqnarray*}
\varepsilon_6&=&\langle 19,29,31\rangle\mid_{76},\\
\varepsilon_7&=&\langle 19,30,31\rangle\mid_{76},
\end{eqnarray*}
both of genus $63$.

Shalom Eliahou and Jean Fromentin found in \cite{EliahouFromentin} numerical semigroups with arbitrarily small Eliahou constant.
In particular, they proved \cite[Proposition 3.1.]{EliahouFromentin} that for positive integers $m,a,b$ with $(3m+1)/2\leq a < b \leq (5m-1)/3$ and for which the elements $a,b,2a,a+b,2b,3a,2a+b,a+2b,3b$ are all different modulo $m$, it holds that 
$EF(m,a,b):=\langle m,a,b\rangle\mid_{4m}$ has Eliahou constant equal to $-1$.
Manuel Delgado constructed in \cite{Delgado}, for each integer number, infinite families of numerical semigroups having Eliahou constant equal to that number.
In particular, he constructed infinite families of semigroups with negative Eliahou constant. 
The semigroups constructed by Delgado are 
$$D^{(i, j)} ( p, \tau ) = \langle m^{(i, j)} , g^{(i, j)} , g^{(i, j)} + 1\rangle \mid_{c^{(i, j)}},$$
for $p$ an even positive integer and $\tau,i,j$ non-negative integers, where
\begin{eqnarray*}
  m^{(i, j)} 
  &=& \frac{p^2}{4}+p(\frac{\tau}{2}+2)+2+j\frac{p}{2} \\
  g^{(i, j)} &=& 
  \frac{p^2}{2}+p(\tau+\frac{7}{2})-\tau  + j(p-1)+ i m^{(i, j)}\\
    c^{(i, j)} &=&
  \frac{p^3}{4}+p^2(\frac{\tau}{2}+2)+2p -\tau + j\frac{p^2}{2}+i\left(\frac{p}{2}+1\right)m^{(i,j)}
\end{eqnarray*}
The Eliahou semigroups in these families found by Eliahou and Formentin and by Delgado still satisfy the Wilf conjecture \cite{EliahouFromentin,Delgado}.
It holds that
\begin{eqnarray*}\varepsilon_1&=&EF(14,22,23)=D^{(0,0)}(4,0)\\
  \varepsilon_2&=&EF(16,25,26)=D^{(0,1)}(4,0)\\
  \varepsilon_3&=&EF(17,26,28)\\
  \varepsilon_4&=&EF(17,27,28)\\
  \varepsilon_5&=&EF(18,28,29)=D^{(0,2)}(4,0)\\
  \varepsilon_6&=&EF(19,29,31)\\
  \varepsilon_7&=&EF(19,30,31)\\
  \end{eqnarray*}
Consequently, the Wilf conjecture holds for all semigroups of genus up to $65$.

Now suppose that $\Lambda$ is not ordinary and let 
$\Lambda'$ be a child of $\Lambda$ obtained by taking away $c+s$ with $s\geq 0$.
Let $p',r',k',q',m',\rho'$ be the parameters involved in the computation of the Eliahou constant for the semigroup $\Lambda'$.
If we define
\begin{eqnarray*}
  \delta_w&=&\left\{\begin{array}{ll}1&\mbox{ if }c+s\mbox{ is a weak generator of }\Lambda\\0&\mbox{ otherwise}\end{array}\right.\\
  \delta_{\rho}&=&\left\{\begin{array}{ll}1&\mbox{ if }\rho\leq s\\0&\mbox{ otherwise,}\end{array}\right.\\
\end{eqnarray*}
then
$m'=m$,
$p'=p-\delta_w$,
$k'=k+s$,
$q'=q+\delta_{\rho}$,
and $\rho'=\rho-s-1+\delta_{\rho}m$.
Let the {\em preceding sibling} $\Lambda'_0$ of $\Lambda'$ be defined as follows. If $c+s$ is the smallest right generator of $\Lambda$, then $\Lambda'_0=\Lambda$. Otherwise, let $s_0$ be such that $c+s_0$ is the largest right generator of $\Lambda$ smaller than $c+s$. Then $\Lambda'_0=\Lambda\setminus\{c+s_0\}$.
Let now
$r'_0$ be the number of right generators of the semigroup $\Lambda'_0$. Then the number of right generators of $\Lambda'$ is $r'=r'_0-\delta_w$.
This enables us to use our exploration of semigroups up to genus $g$ for finding all Eliahou semigroups of genus up to $g+1$.
In fact, our algorithm only checks the Eliahou inequality for semigroups with $k\geq 3$.
But as proved in \cite[Section 7.2]{eliahou} this covers all possible counterexamples.

With the new algorithm we found three further Eliahou semigroups of genus $67$.
They are
\begin{eqnarray*}
\varepsilon_8&=&\langle 20, 31, 32\rangle\mid_{80}\\
\varepsilon_9&=&\langle 20, 32, 33\rangle\mid_{80}\\
\varepsilon_{10}&=&\langle 19, 26, 27\rangle\mid_{90}\\
\end{eqnarray*}
The three semigroups $\varepsilon_8$, $\varepsilon_9$, $\varepsilon_{10}$,
have Eliahou constant equal to $-1$.
The parameters involved in the computation of the Eliahou constant
of the semigroups $\varepsilon_8$ and $\varepsilon_9$
coincide and they are
$$c=80,\ m=20,\ q=4,\ \rho=0,\ p=13,\ r=10,\ k=13.$$
The parameters involved in the computation of the Eliahou constant
of the semigroup $\varepsilon_{10}$ are
$$c=90,\ m=19,\ q=5,\ \rho=5,\ p=7,\ r=4,\ k=23.$$
The semigroups $\varepsilon_8$ and $\varepsilon_9$ are of the type described by Eliahou and Fromentin in \cite{EliahouFromentin}. 
Furthermore, the semigroup $\varepsilon_8$ is the Delgado semigroup $D^{(0,3)}(4,0)$. That is,
\begin{eqnarray*}\varepsilon_8&=&EF(20,31,32)=D^{(0,3)}(4,0)\\
  \varepsilon_9&=&EF(20,32,33)\\
\end{eqnarray*}
The semigroup $\varepsilon_{10}$ is neither of Eliahou-Fromentin type nor of Delgado's type.

Shalom Eliahou noticed that, indeed, $\varepsilon_{10}$ is a member of a new family of Eliahou semigroups, defined as
$$BEF_t=\langle 2t+1, 3t-1, 3t \rangle\mid_{10t},$$
for $t\geq 9$.
All the semigroups in this family satisfy
$$c=10t,\ m=2t+1,\ q=5,\ \rho=5,$$
and, since
\begin{eqnarray*}
{BEF}_t&=&\{0,2t+1,3t-1,3t,\    
4t+2,5t,5t+1,6t-2,6t-1,6t,\ 
6t+3,\                      
\\&&7t+1,7t+2,\                 
8t-1, 8t, 8t+1,\            
8t+4,\                      
9t-3, 9t-2, 9t-1, 9t,\      
9t+2,9t+3,\                 
10t,\dots\
\},\end{eqnarray*}
the rank is $k=23$ whenever $t\geq 8$. In particular, their genus is $10t-23$.
The elements between the conductor $c$ and $c+m-1$ that are not generators
are
$$\
10t, 10t+1, 10t+2,\                
10t+5,             \               
11t-2, 11t-1, 11t, 11t+1,\         
11t+3, 11t+4,             \        
12t-4, 12t-3, 12t-2, 12t-1, 12t.\  
\                                  
\                                  
$$
For $t\geq 9$ this list contains exactly $15$ elements, and so the number of right generators is $r=2t+1-15=2t-14$.
For $t=8$, the list contains $14$ elements and so $r=2t-13$.
Hence, $E(BEF_8) =4$ while $E(BEF_t) = -1$ for all $t\geq 9$.
This shows that the family $BEF_t$ for $t\geq 9$ is indeed a family of Eliahou semigroups and that the condition $t\geq 9$ is necessary.
It is easy to check that the Wilf conjecture holds for these semigroups since $c=10t$, $k=23$, $p=2t-11$.
One can check that \begin{eqnarray*}\varepsilon_{10}&=&BEF_9.\end{eqnarray*}

We can not ensure a priori that the three Eliahou semigroups $\varepsilon_8,\varepsilon_9,\varepsilon_{10}$ are {\em all} Eliahou semigroups of genus between $66$ and $67$ because we have the integer bit length limitation. We will analyze it in next section.

    \mut{
\begin{itemize}
\item
\begin{alltt}
Lambda={ 0 19 26 27 38 45 46 52 53 54 57 64 65 71 72 73 76 78 79 80 81 83 84 90 ... }
Is semigroup?
Yes!
Generators?
19 26 27 93 94 96 101 
c=90 m=19 u=7 v=1 g=67
p=7 r=4 l=23 q=5 rho=5
Wilf inequality: 161 >= 90
Eliahou inequality: -1 >= 0
\end{alltt}

\item
\begin{alltt}
Lambda={ 0 20 31 32 40 51 52 60 62 63 64 71 72 80 ... }
Is semigroup?
Yes!
Generators?
20 31 32 81 85 86 87 88 89 90 97 98 99 
c=80 m=20 u=11 v=1 g=67
p=13 r=10 l=13 q=4 rho=0
Wilf inequality: 169 >= 80
Eliahou inequality: -1 >= 0
\end{alltt}
\item
\begin{alltt}
Lambda={ 0 20 32 33 40 52 53 60 64 65 66 72 73 80 ... }
Is semigroup?
Yes!
Generators?
20 32 33 81 82 83 87 88 89 90 91 94 95 
c=80 m=20 u=12 v=1 g=67
p=13 r=10 l=13 q=4 rho=0
Wilf inequality: 169 >= 80
Eliahou inequality: -1 >= 0
\end{alltt}

\end{itemize}
    }

\section{How far can we get using $128$ bit integers?}
\label{s:howfar}

Let $n_g$ be the number of semigroups of genus $g$.
Let $E_g$ be the set of Eliahou counterexamples of genus smaller than or equal to $g$.
As a consequence of the previous sections, if we obtain the parameters $G$ and $S$ of all semigroups of genus up to $g$, then we can also deduce $n_{g+3}$ and $E_{g+1}$.
The limitation of our algorithm comes from the fact that the variables $G$ and $S$ allocate values in positions $0$ to $c-1$.
If we use $128$ bit integers, this means that we explore correctly the semigroups of conductor up to $128$.

\subsection{Exploring up to genus $64$}

Since $c\leq 2g$, if we use $128$ bit integers, this guarantees that the algorithm will proceed correctly to explore the tree up to genus $64$, thus guaranteeing the correct computation of $n_{67}$ and $E_{65}$.

\subsection{Exploring up to genus $65$}

If we use the algorithm to explore up to genus $65$, we will missexplore semigroups of genus $65$ and conductor $129$ or $130$, that is, we will be missexploring exactly the symmetric (i.e. rank equal to $g$) and pseudo-symmetric (i.e. rank equal to $g-1$) semigroups of genus $65$.
The error in these cases will be that we will not be able to write
\begin{itemize}
\item $S_{c-1}=1$ when $c\geq 129$
\item $S_{c-2}=1$ and $G_{c-2}=1$ when $c\geq 130$
\end{itemize}
But the parameters such as the genus, conductor, multiplicity, first jump and second jump, or the number of right generators will be correctly computed.

\subsubsection*{Computation of $n_{68}$}

To analyze this case we need the next lemmas which are consequences of results one can find in \cite{BrBu}.

\begin{lemma}\cite[Lemma 3]{BrBu}
Non-hyperelliptic symmetric semigroups have no children.
\end{lemma}

\begin{lemma}
Pseudo-symmetric semigroups have no grandchildren.
\end{lemma}
  \begin{proof}
From \cite[Lemma 4]{BrBu} we deduce that there is a unique pseudo-semigroup of genus $g$ with only one interval of
non-gaps between $0$ and the conductor. If $g\geq 5$ then this semigroup has only one child and no grandchildren.
From \cite[Lemma 5]{BrBu} we deduce that pseudo-semigroups with $\lambda_1=3$ have one child and no grandchildren.
Finally, \cite[Lemma 6]{BrBu} says that each pseudo-symmetric semigroup with $\lambda_1\neq 3$ and with more than one interval of non-gaps between 0 and the conductor is a leaf in the semigroup tree.
  \end{proof}

From these lemmas we deduce that the missexplored symmetric and pseudo-symmetric semigroups of genus $65$ do not interfere in the computation of $n_{68}$.
Indeed, these semigroups have at most one child and no grandchildren.
Hence they do not contribute to $n_{68}$.

\subsubsection*{Computation of $E_{66}$}

As a consequence of \cite[Lemma 4]{BrBu}, \cite[Lemma 5]{BrBu}, and \cite[Lemma 6]{BrBu}, we can state the following lemma.  
  \begin{lemma}
    The unique semigroups that are children of symmetric or pseudo-symmetric semigroups of a given genus $g\geq 5$ are:
    \begin{enumerate}
    \item The hyperelliptic semigroup of genus $g+1$
    \item The semigroup $\{0,g,g+1,\dots,2g-3\}\cup [2g,\infty)$
    \item The semigroup $\{0,3,6,\dots,3t,3(t+1)-1,3(t+1),3(t+2)-1,3(t+2),\dots, 3(2t-1)-1,3(2t-1)\}\cup\{3(2t-1)+2,3(2t-1)+3\}\cup [3(2t-1)+5,\infty)$ for an adequate $t$
    \item The semigroup $\{0,3,6,\dots,3t,3(t+1),3(t+1)+1,3(t+2),3(t+2)+1,\dots, 3(2t),3(2t)+1\}\cup\{3(2t)+3,3(2t)+4\}\cup[3(2t)+6,\infty)$ for an adequate $t$.
      \end{enumerate}
  \end{lemma}

  Now we can check that all semigroups in the previous lemma have non-negative Eliahou constant, thus proving that our algorithm using 128-bit integers guarantees the computation of $E_{66}$.
  Indeed,
  \begin{enumerate}
  \item The hyperelliptic semigroup of genus $g+1$ has parameters
          $m=2$, $p=2$, $k=g+1$, $r=1$, $q=g+1$, and $\rho=0$.
    Hence, the Eliahou constant is
    $(g+1)-(g+1)=0$.
  \item The semigroup $\{0,g,g+1,\dots,2g-3\}\cup [2g,\infty)$
    has parameters
    $m=g$, $p=g-2$, $k=g-1$, $r=0$, $q=2$, and $\rho=0$.
    Hence, the Eliahou constant is
    $k(p-r)-q(m-r)+\rho=(g-1)(g-2)-2g=g^2-5g+2$, which is positive for $g\geq 5$.    
  \item The semigroup $\{0,3,6,\dots,3t,3(t+1)-1,3(t+1),3(t+2)-1,3(t+2),\dots, 3(2t-1)-1,3(2t-1)\}\cup\{3(2t-1)+2,3(2t-1)+3\}\cup [3(2t-1)+5,\infty)$ for an adequate $t$
    has parameters
        $m=3$, $p=2$, $k=3t+1$, $r=0$, $q=2t+1$, and $\rho=1$.
    Hence, the Eliahou constant is
    $k(p-r)-q(m-r)+\rho=(3t+1)2-(2t+1)3+1=0$.
  \item The semigroup $\{0,3,6,\dots,3t,3(t+1),3(t+1)+1,3(t+2),3(t+2)+1,\dots, 3(2t),3(2t)+1\}\cup\{3(2t)+3,3(2t)+4\}\cup[3(2t)+6,\infty)$ for an adequate $t$.
    has parameters
    $m=3$, $p=2$, $k=3t+3$, $r=0$, $q=2t+2$, and $\rho=0$.
    Hence, the Eliahou constant is
    $k(p-r)-q(m-r)+\rho=(3t+3)2-(2t+2)3=0$.
  \end{enumerate}

  These results allow us to conclude with the next lemmas.

  \begin{proposition}
    The unique Eliahou semigroups of genus up to $66$ are exactly $\varepsilon_1,\varepsilon_2,\varepsilon_3,\varepsilon_4,\varepsilon_5,\varepsilon_6,\varepsilon_7$.   
    \end{proposition}
  
  \begin{corollary}
The Wilf conjecture holds for all semigroups of genus up to 66.
    \end{corollary}
  
\subsection{Exploring up to genus $66$}

At this point we can say that we found all semigroups of $E_{67}$ except the elements of $E_{67}$ that have genus $67$ and are children of semigroups of genus $66$ of rank $g,g-1,g-2,g-3$, should they exist. The children of rank-$g$ and rank-$(g-1)$ semigroups can be discarded from the results in the previous section. It remains as an open question to see whether there may be children of rank-$(g-2)$ and rank-$(g-3)$ semigroups with negative Eliahou constant.

\section*{Acknowledgment}
The author would like to thank Julio Fernández-González for many helpful insights and comments and César Marín Rodríguez for his ideas and his contribution to the implementation of the algorithm.
She would also like to thank Manuel Delgado and Shalom Eliahou for their comments and for noticing that two of the Eliahou semigroups of genus $67$ are, in fact, Eliahou-Fromentin semigroups.
Also, as mentioned above, Shalom Eliahou observed that a new family of semigroups could be defined with Eliahou constant equal to $-1$, containing the semigroup
$\langle 19, 26, 27\rangle\mid_{90}$.
The author would also like to thank Paul Shin for his careful reading of the manuscript and for his interesting observations as well as the referee that helped improving the proof of Lemma 3.

All the graphs have been drawn using the {\tt drawsgtree} tool, which can be downloaded from {\tt https://github.com/mbrasamoros/drawsgtree}.

The author was supported by the Spanish government under grant PID2021-124928NB-I00, and by the Catalan government under grant 2021 SGR 00115.

\bibliographystyle{plain}
\bibliography{bib}

\begin{thebibliography}{10}

\bibitem{Br:fibonacci}
M.~Bras-Amor{\'o}s.
\newblock Fibonacci-like behavior of the number of numerical semigroups of a
  given genus.
\newblock {\em Semigroup Forum}, 76(2):379--384, 2008.

\bibitem{Br:bounds}
M.~Bras-Amor{\'o}s.
\newblock Bounds on the number of numerical semigroups of a given genus.
\newblock {\em J. Pure Appl. Algebra}, 213(6):997--1001, 2009.

\bibitem{BrBu}
M.~Bras-Amor{\'o}s and S.~Bulygin.
\newblock Towards a better understanding of the semigroup tree.
\newblock {\em Semigroup Forum}, 79(3):561--574, 2009.

\bibitem{RGD-code}
M.~Bras-Amor\'os and J.~Fern\'andez-Gonz\'alez.
\newblock {\footnotesize\texttt{https://github.com/mbrasamoros/RGD-algorithm}}.

\bibitem{seeds}
M.~Bras-Amor{\'o}s and J.~Fern{\'a}ndez-Gonz{\'a}lez.
\newblock Computation of numerical semigroups by means of seeds.
\newblock {\em Math. Comp.}, 87(313):2539--2550, 2018.

\bibitem{rgd}
M.~Bras-Amor\'{o}s and J.~Fern\'{a}ndez-Gonz\'{a}lez.
\newblock The right-generators descendant of a numerical semigroup.
\newblock {\em Math. Comp.}, 89(324):2017--2030, 2020.

\bibitem{mdai}
M.~Bras-Amorós and C.~Marín~Rodríguez.
\newblock New {E}liahou semigroups and verification of the {W}ilf conjecture
  for genus up to 65.
\newblock In Torra and Narukawa, editors, {\em Modeling Decisions for
  Artificial Intelligence}, volume 12898 of {\em Lect. Notes Comput. Sci.},
  pages 17--27. Springer, 2021.

\bibitem{Delgado}
M.~Delgado.
\newblock On a question of {E}liahou and a conjecture of {W}ilf.
\newblock {\em Math.Z.}, 288(1-2):595--627, 2018.

\bibitem{DelgadoTrimming}
M.~Delgado.
\newblock Trimming the numerical semigroups tree to probe {W}ilf's conjecture
  to higher genus, 2019.

\bibitem{DelgadoIMNS}
M.~Delgado.
\newblock {\em Conjecture of Wilf: a survey}, volume~40 of {\em Springer INdAM
  Series}, pages 39--62.
\newblock Springer International Publishing, 2020.

\bibitem{eliahou}
S.~Eliahou.
\newblock Wilf's conjecture and {M}acaulay's theorem.
\newblock {\em J. Eur. Math. Soc. (JEMS)}, 20(9):2105--2129, 2018.

\bibitem{EliahouFromentin}
S.~Eliahou and J.~Fromentin.
\newblock Near-misses in {W}ilf's conjecture.
\newblock {\em Semigroup Forum}, 98(2):285--298, 2019.

\bibitem{FH-code}
J.~Fromentin and F.~Hivert.
\newblock {\footnotesize\texttt{https://github.com/hivert/NumericMonoid}}.

\bibitem{FromentinHivert}
J.~Fromentin and F.~Hivert.
\newblock Exploring the tree of numerical semigroups.
\newblock {\em Math. Comp.}, 85(301):2553--2568, 2016.

\bibitem{RGS}
J.~C. Rosales and P.~A. Garc{\'{\i}}a-S{\'a}nchez.
\newblock {\em Numerical semigroups}, volume~20 of {\em Developments in
  Mathematics}.
\newblock Springer, New York, 2009.

\bibitem{Wilf}
H.~S. Wilf.
\newblock A circle-of-lights algorithm for the ``money-changing problem''.
\newblock {\em Amer. Math. Monthly}, 85(7):562--565, 1978.

\end{thebibliography}

\newpage
\appendix

\section{Illustrative examples}

The examples in this appendix are included in order to illustrate all cases considered in Lemma~\ref{l:ggc}.

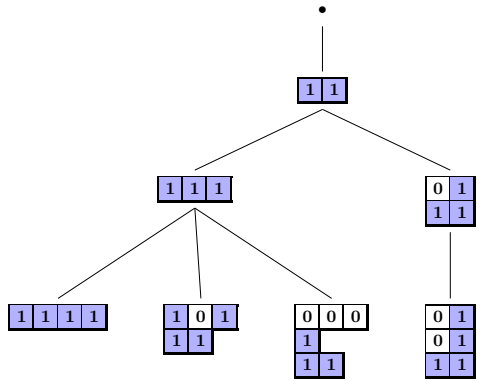
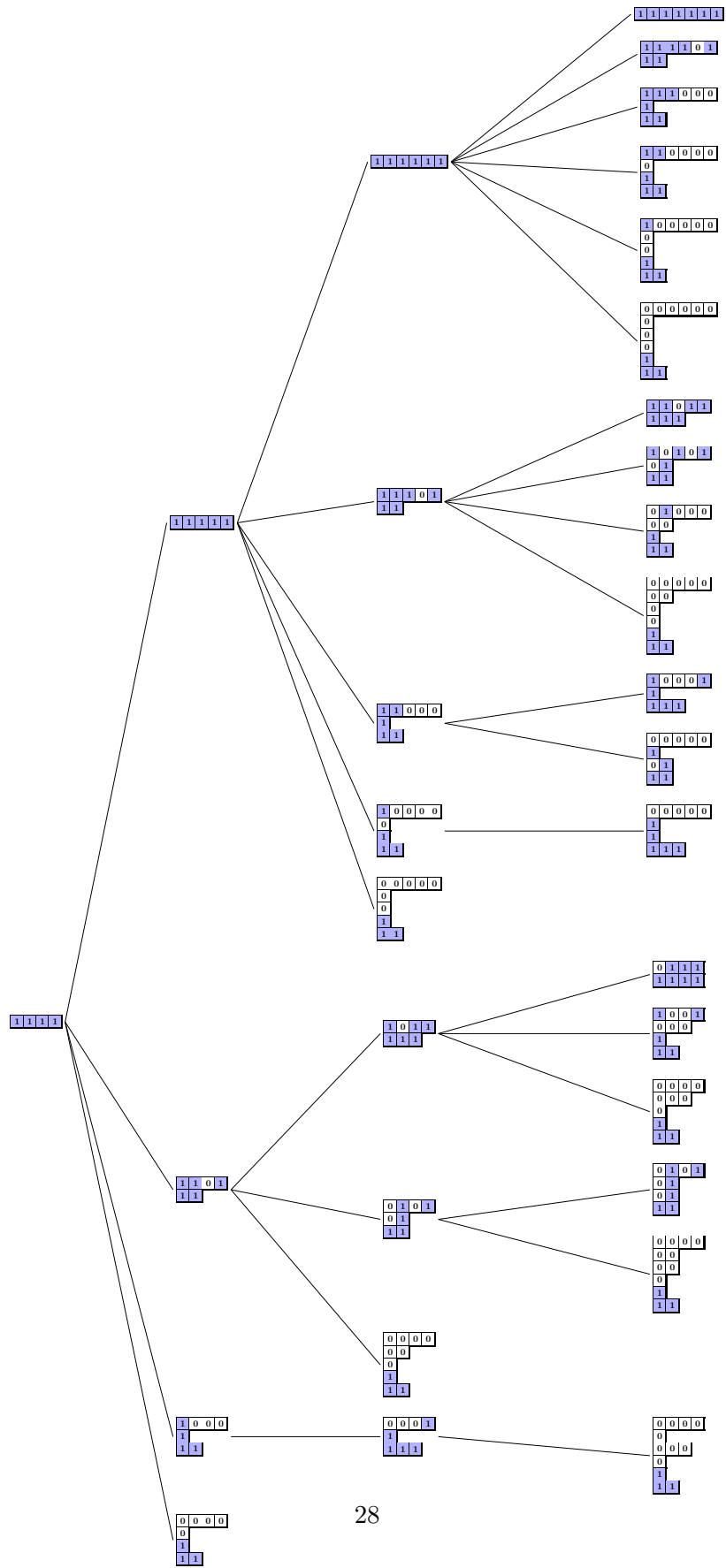
\begin{figure}[h]
  \begin{center}
%
%
%
%

\newcolumntype{M}{>{\centering\arraybackslash}m{.5cm}}\setlength\tabcolsep{0pt}\setlength\arrayrulewidth{1pt}\providecommand\circledcolorednumb{}\renewcommand\circledcolorednumb[2]{\resizebox{0.066975\textwidth}{!}{\tikz[baseline=(char.center)]{\node[shape = circle,draw, inner sep = 2pt,fill=#1](char)    {\phantom{00}};\node[anchor=center] at (char.center) {\makebox(0,0){\large{#2}}};}}}
\robustify{\circledcolorednumb}
\providecommand\nongap{}\renewcommand\nongap[1]{\circledcolorednumb{yellow}{{\bf#1}}}
\providecommand\gap{}\renewcommand\gap[1]{\circledcolorednumb{white}{\phantom{#1}}}
\providecommand\generator{}\renewcommand\generator[1]{\circledcolorednumb{orange}{{\bf#1}}}
\providecommand\seed{}\renewcommand\seed[1]{\circledcolorednumb{blue!30}{{\bf#1}}}
\providecommand\nonseed{}\renewcommand\nonseed[1]{\circledcolorednumb{white}{{\bf#1}}}
\providecommand\dotscircles{}\renewcommand\dotscircles{\resizebox{0.066975\textwidth}{!}{\dots}}
\providecommand\gapingapset{}\renewcommand\gapingapset[1]{\circledcolorednumb{green!30}{#1}}
\providecommand\nongapingapset{}\renewcommand\nongapingapset[1]{\phantom{\gapingapset{#1}}}
\providecommand\coloredseed{}\renewcommand\coloredseed{\cellcolor{blue!30}}
\scalebox{0.600000}{\adjustbox{max width=\textwidth,max height=.9\textheight}{\begin{tikzpicture}[grow=down,sibling distance=10.000000mm]\tikzset{every tree node/.style={anchor=north}}\tikzset{level 1/.style={level distance=1.750000cm}}\tikzset{level 2/.style={level distance=2.187500cm}}\tikzset{level 3/.style={level distance=2.843750cm}}\tikzset{level 4/.style={level distance=3.500000cm}}\tikzset{level 5/.style={level distance=4.375000cm}}\tikzset{level 6/.style={level distance=4.462500cm}}\tikzset{level 7+/.style={level distance=4.725000cm}}\Tree[.{\begin{tabular}{c}{$\bullet$} \\\end{tabular}} [.{\begin{tabular}{c}{{\bf \begin{tabular}{|@{\rule[-.15cm]{0pt}{.5cm}}*{2}{M |}}\hhline{|-|-|}
  \coloredseed 1  & \coloredseed 1  \\\hhline{|-|-|}
\end{tabular}}} \\\end{tabular}} [.{\begin{tabular}{c}{{\bf \begin{tabular}{|@{\rule[-.15cm]{0pt}{.5cm}}*{3}{M |}}\hhline{|-|-|-|}
  \coloredseed 1  & \coloredseed 1  & \coloredseed 1  \\\hhline{|-|-|-|}
\end{tabular}}} \\\end{tabular}} [.{\begin{tabular}{c}{{\bf \begin{tabular}{|@{\rule[-.15cm]{0pt}{.5cm}}*{4}{M |}}\hhline{|-|-|-|-|}
  \coloredseed 1  & \coloredseed 1  & \coloredseed 1  & \coloredseed 1  \\\hhline{|-|-|-|-|}
\end{tabular}}} \\\end{tabular}} ][.{\begin{tabular}{c}{{\bf \begin{tabular}{|@{\rule[-.15cm]{0pt}{.5cm}}*{3}{M |}}\hhline{|-|-|-|}
  \coloredseed 1  & \cellcolor{white} 0  & \coloredseed 1  \\\hhline{|-|-|-|}
  \coloredseed 1  & \coloredseed 1  \\\hhline{|-|-|}
\end{tabular}}} \\\end{tabular}} ][.{\begin{tabular}{c}{{\bf \begin{tabular}{|@{\rule[-.15cm]{0pt}{.5cm}}*{3}{M |}}\hhline{|-|-|-|}
  \cellcolor{white} 0  & \cellcolor{white} 0  & \cellcolor{white} 0  \\\hhline{|-|-|-|}
  \coloredseed 1  \\\hhline{|-|-|}
  \coloredseed 1  & \coloredseed 1  \\\hhline{|-|-|}
\end{tabular}}} \\\end{tabular}} ]][.{\begin{tabular}{c}{{\bf \begin{tabular}{|@{\rule[-.15cm]{0pt}{.5cm}}*{2}{M |}}\hhline{|-|-|}
  \cellcolor{white} 0  & \coloredseed 1  \\\hhline{|-|-|}
  \coloredseed 1  & \coloredseed 1  \\\hhline{|-|-|}
\end{tabular}}} \\\end{tabular}} [.{\begin{tabular}{c}{{\bf \begin{tabular}{|@{\rule[-.15cm]{0pt}{.5cm}}*{2}{M |}}\hhline{|-|-|}
  \cellcolor{white} 0  & \coloredseed 1  \\\hhline{|-|-|}
  \cellcolor{white} 0  & \coloredseed 1  \\\hhline{|-|-|}
  \coloredseed 1  & \coloredseed 1  \\\hhline{|-|-|}
\end{tabular}}} \\\end{tabular}} ]]]]\end{tikzpicture}}}
  \end{center}
  \caption{Seed tables of the children, grandchildren, and great grandchildren of  $\{0, 1, \dots\}$}\label{fig:exa}
\end{figure}

\begin{figure}[h]
  \begin{center}
%
%
%
%

\newcolumntype{M}{>{\centering\arraybackslash}m{.5cm}}\setlength\tabcolsep{0pt}\setlength\arrayrulewidth{1pt}\providecommand\circledcolorednumb{}\renewcommand\circledcolorednumb[2]{\resizebox{0.100463\textwidth}{!}{\tikz[baseline=(char.center)]{\node[shape = circle,draw, inner sep = 2pt,fill=#1](char)    {\phantom{00}};\node[anchor=center] at (char.center) {\makebox(0,0){\large{#2}}};}}}
\robustify{\circledcolorednumb}
\providecommand\nongap{}\renewcommand\nongap[1]{\circledcolorednumb{yellow}{{\bf#1}}}
\providecommand\gap{}\renewcommand\gap[1]{\circledcolorednumb{white}{\phantom{#1}}}
\providecommand\generator{}\renewcommand\generator[1]{\circledcolorednumb{orange}{{\bf#1}}}
\providecommand\seed{}\renewcommand\seed[1]{\circledcolorednumb{blue!30}{{\bf#1}}}
\providecommand\nonseed{}\renewcommand\nonseed[1]{\circledcolorednumb{white}{{\bf#1}}}
\providecommand\dotscircles{}\renewcommand\dotscircles{\resizebox{0.100463\textwidth}{!}{\dots}}
\providecommand\gapingapset{}\renewcommand\gapingapset[1]{\circledcolorednumb{green!30}{#1}}
\providecommand\nongapingapset{}\renewcommand\nongapingapset[1]{\phantom{\gapingapset{#1}}}
\providecommand\coloredseed{}\renewcommand\coloredseed{\cellcolor{blue!30}}
\scalebox{0.700000}{\adjustbox{max width=\textwidth,max height=.9\textheight}{\begin{tikzpicture}[grow=down,sibling distance=5.000000mm]\tikzset{every tree node/.style={anchor=north}}\tikzset{level 1/.style={level distance=2.100000cm}}\tikzset{level 2/.style={level distance=2.625000cm}}\tikzset{level 3/.style={level distance=3.412500cm}}\tikzset{level 4/.style={level distance=4.200000cm}}\tikzset{level 5/.style={level distance=5.250000cm}}\tikzset{level 6/.style={level distance=5.355000cm}}\tikzset{level 7+/.style={level distance=5.670000cm}}\Tree[.{\begin{tabular}{c}{{\bf \begin{tabular}{|@{\rule[-.15cm]{0pt}{.5cm}}*{2}{M |}}\hhline{|-|-|}
  \coloredseed 1  & \coloredseed 1  \\\hhline{|-|-|}
\end{tabular}}} \\\end{tabular}} [.{\begin{tabular}{c}{{\bf \begin{tabular}{|@{\rule[-.15cm]{0pt}{.5cm}}*{3}{M |}}\hhline{|-|-|-|}
  \coloredseed 1  & \coloredseed 1  & \coloredseed 1  \\\hhline{|-|-|-|}
\end{tabular}}} \\\end{tabular}} [.{\begin{tabular}{c}{{\bf \begin{tabular}{|@{\rule[-.15cm]{0pt}{.5cm}}*{4}{M |}}\hhline{|-|-|-|-|}
  \coloredseed 1  & \coloredseed 1  & \coloredseed 1  & \coloredseed 1  \\\hhline{|-|-|-|-|}
\end{tabular}}} \\\end{tabular}} [.{\begin{tabular}{c}{{\bf \begin{tabular}{|@{\rule[-.15cm]{0pt}{.5cm}}*{5}{M |}}\hhline{|-|-|-|-|-|}
  \coloredseed 1  & \coloredseed 1  & \coloredseed 1  & \coloredseed 1  & \coloredseed 1  \\\hhline{|-|-|-|-|-|}
\end{tabular}}} \\\end{tabular}} ][.{\begin{tabular}{c}{{\bf \begin{tabular}{|@{\rule[-.15cm]{0pt}{.5cm}}*{4}{M |}}\hhline{|-|-|-|-|}
  \coloredseed 1  & \coloredseed 1  & \cellcolor{white} 0  & \coloredseed 1  \\\hhline{|-|-|-|-|}
  \coloredseed 1  & \coloredseed 1  \\\hhline{|-|-|}
\end{tabular}}} \\\end{tabular}} ][.{\begin{tabular}{c}{{\bf \begin{tabular}{|@{\rule[-.15cm]{0pt}{.5cm}}*{4}{M |}}\hhline{|-|-|-|-|}
  \coloredseed 1  & \cellcolor{white} 0  & \cellcolor{white} 0  & \cellcolor{white} 0  \\\hhline{|-|-|-|-|}
  \coloredseed 1  \\\hhline{|-|-|}
  \coloredseed 1  & \coloredseed 1  \\\hhline{|-|-|}
\end{tabular}}} \\\end{tabular}} ][.{\begin{tabular}{c}{{\bf \begin{tabular}{|@{\rule[-.15cm]{0pt}{.5cm}}*{4}{M |}}\hhline{|-|-|-|-|}
  \cellcolor{white} 0  & \cellcolor{white} 0  & \cellcolor{white} 0  & \cellcolor{white} 0  \\\hhline{|-|-|-|-|}
  \cellcolor{white} 0  \\\hhline{|-|}
  \coloredseed 1  \\\hhline{|-|-|}
  \coloredseed 1  & \coloredseed 1  \\\hhline{|-|-|}
\end{tabular}}} \\\end{tabular}} ]][.{\begin{tabular}{c}{{\bf \begin{tabular}{|@{\rule[-.15cm]{0pt}{.5cm}}*{3}{M |}}\hhline{|-|-|-|}
  \coloredseed 1  & \cellcolor{white} 0  & \coloredseed 1  \\\hhline{|-|-|-|}
  \coloredseed 1  & \coloredseed 1  \\\hhline{|-|-|}
\end{tabular}}} \\\end{tabular}} [.{\begin{tabular}{c}{{\bf \begin{tabular}{|@{\rule[-.15cm]{0pt}{.5cm}}*{3}{M |}}\hhline{|-|-|-|}
  \cellcolor{white} 0  & \coloredseed 1  & \coloredseed 1  \\\hhline{|-|-|-|}
  \coloredseed 1  & \coloredseed 1  & \coloredseed 1  \\\hhline{|-|-|-|}
\end{tabular}}} \\\end{tabular}} ][.{\begin{tabular}{c}{{\bf \begin{tabular}{|@{\rule[-.15cm]{0pt}{.5cm}}*{3}{M |}}\hhline{|-|-|-|}
  \cellcolor{white} 0  & \cellcolor{white} 0  & \cellcolor{white} 0  \\\hhline{|-|-|-|}
  \cellcolor{white} 0  & \cellcolor{white} 0  \\\hhline{|-|-|}
  \coloredseed 1  \\\hhline{|-|-|}
  \coloredseed 1  & \coloredseed 1  \\\hhline{|-|-|}
\end{tabular}}} \\\end{tabular}} ]][.{\begin{tabular}{c}{{\bf \begin{tabular}{|@{\rule[-.15cm]{0pt}{.5cm}}*{3}{M |}}\hhline{|-|-|-|}
  \cellcolor{white} 0  & \cellcolor{white} 0  & \cellcolor{white} 0  \\\hhline{|-|-|-|}
  \coloredseed 1  \\\hhline{|-|-|}
  \coloredseed 1  & \coloredseed 1  \\\hhline{|-|-|}
\end{tabular}}} \\\end{tabular}} ]][.{\begin{tabular}{c}{{\bf \begin{tabular}{|@{\rule[-.15cm]{0pt}{.5cm}}*{2}{M |}}\hhline{|-|-|}
  \cellcolor{white} 0  & \coloredseed 1  \\\hhline{|-|-|}
  \coloredseed 1  & \coloredseed 1  \\\hhline{|-|-|}
\end{tabular}}} \\\end{tabular}} [.{\begin{tabular}{c}{{\bf \begin{tabular}{|@{\rule[-.15cm]{0pt}{.5cm}}*{2}{M |}}\hhline{|-|-|}
  \cellcolor{white} 0  & \coloredseed 1  \\\hhline{|-|-|}
  \cellcolor{white} 0  & \coloredseed 1  \\\hhline{|-|-|}
  \coloredseed 1  & \coloredseed 1  \\\hhline{|-|-|}
\end{tabular}}} \\\end{tabular}} [.{\begin{tabular}{c}{{\bf \begin{tabular}{|@{\rule[-.15cm]{0pt}{.5cm}}*{2}{M |}}\hhline{|-|-|}
  \cellcolor{white} 0  & \coloredseed 1  \\\hhline{|-|-|}
  \cellcolor{white} 0  & \coloredseed 1  \\\hhline{|-|-|}
  \cellcolor{white} 0  & \coloredseed 1  \\\hhline{|-|-|}
  \coloredseed 1  & \coloredseed 1  \\\hhline{|-|-|}
\end{tabular}}} \\\end{tabular}} ]]]]\end{tikzpicture}}}
  \end{center}
  \caption{Seed tables of the children, grandchildren, and great grandchildren of  $\{0, 2, \dots\}$}\label{fig:exb}
\end{figure}

\begin{figure}[h]
  \begin{center}
%
%
%
%

\newcolumntype{M}{>{\centering\arraybackslash}m{.5cm}}\setlength\tabcolsep{0pt}\setlength\arrayrulewidth{1pt}\providecommand\circledcolorednumb{}\renewcommand\circledcolorednumb[2]{\resizebox{0.133951\textwidth}{!}{\tikz[baseline=(char.center)]{\node[shape = circle,draw, inner sep = 2pt,fill=#1](char)    {\phantom{00}};\node[anchor=center] at (char.center) {\makebox(0,0){\large{#2}}};}}}
\robustify{\circledcolorednumb}
\providecommand\nongap{}\renewcommand\nongap[1]{\circledcolorednumb{yellow}{{\bf#1}}}
\providecommand\gap{}\renewcommand\gap[1]{\circledcolorednumb{white}{\phantom{#1}}}
\providecommand\generator{}\renewcommand\generator[1]{\circledcolorednumb{orange}{{\bf#1}}}
\providecommand\seed{}\renewcommand\seed[1]{\circledcolorednumb{blue!30}{{\bf#1}}}
\providecommand\nonseed{}\renewcommand\nonseed[1]{\circledcolorednumb{white}{{\bf#1}}}
\providecommand\dotscircles{}\renewcommand\dotscircles{\resizebox{0.133951\textwidth}{!}{\dots}}
\providecommand\gapingapset{}\renewcommand\gapingapset[1]{\circledcolorednumb{green!30}{#1}}
\providecommand\nongapingapset{}\renewcommand\nongapingapset[1]{\phantom{\gapingapset{#1}}}
\providecommand\coloredseed{}\renewcommand\coloredseed{\cellcolor{blue!30}}
\adjustbox{max width=\textwidth,max height=.9\textheight}{\begin{tikzpicture}[grow=down,sibling distance=3.000000mm]\tikzset{every tree node/.style={anchor=north}}\tikzset{level 1/.style={level distance=2.800000cm}}\tikzset{level 2/.style={level distance=3.500000cm}}\tikzset{level 3/.style={level distance=4.550000cm}}\tikzset{level 4/.style={level distance=5.600000cm}}\tikzset{level 5/.style={level distance=7.000000cm}}\tikzset{level 6/.style={level distance=7.140000cm}}\tikzset{level 7+/.style={level distance=7.560000cm}}\Tree[.{\begin{tabular}{c}{{\bf \begin{tabular}{|@{\rule[-.15cm]{0pt}{.5cm}}*{3}{M |}}\hhline{|-|-|-|}
  \coloredseed 1  & \coloredseed 1  & \coloredseed 1  \\\hhline{|-|-|-|}
\end{tabular}}} \\\end{tabular}} [.{\begin{tabular}{c}{{\bf \begin{tabular}{|@{\rule[-.15cm]{0pt}{.5cm}}*{4}{M |}}\hhline{|-|-|-|-|}
  \coloredseed 1  & \coloredseed 1  & \coloredseed 1  & \coloredseed 1  \\\hhline{|-|-|-|-|}
\end{tabular}}} \\\end{tabular}} [.{\begin{tabular}{c}{{\bf \begin{tabular}{|@{\rule[-.15cm]{0pt}{.5cm}}*{5}{M |}}\hhline{|-|-|-|-|-|}
  \coloredseed 1  & \coloredseed 1  & \coloredseed 1  & \coloredseed 1  & \coloredseed 1  \\\hhline{|-|-|-|-|-|}
\end{tabular}}} \\\end{tabular}} [.{\begin{tabular}{c}{{\bf \begin{tabular}{|@{\rule[-.15cm]{0pt}{.5cm}}*{6}{M |}}\hhline{|-|-|-|-|-|-|}
  \coloredseed 1  & \coloredseed 1  & \coloredseed 1  & \coloredseed 1  & \coloredseed 1  & \coloredseed 1  \\\hhline{|-|-|-|-|-|-|}
\end{tabular}}} \\\end{tabular}} ][.{\begin{tabular}{c}{{\bf \begin{tabular}{|@{\rule[-.15cm]{0pt}{.5cm}}*{5}{M |}}\hhline{|-|-|-|-|-|}
  \coloredseed 1  & \coloredseed 1  & \coloredseed 1  & \cellcolor{white} 0  & \coloredseed 1  \\\hhline{|-|-|-|-|-|}
  \coloredseed 1  & \coloredseed 1  \\\hhline{|-|-|}
\end{tabular}}} \\\end{tabular}} ][.{\begin{tabular}{c}{{\bf \begin{tabular}{|@{\rule[-.15cm]{0pt}{.5cm}}*{5}{M |}}\hhline{|-|-|-|-|-|}
  \coloredseed 1  & \coloredseed 1  & \cellcolor{white} 0  & \cellcolor{white} 0  & \cellcolor{white} 0  \\\hhline{|-|-|-|-|-|}
  \coloredseed 1  \\\hhline{|-|-|}
  \coloredseed 1  & \coloredseed 1  \\\hhline{|-|-|}
\end{tabular}}} \\\end{tabular}} ][.{\begin{tabular}{c}{{\bf \begin{tabular}{|@{\rule[-.15cm]{0pt}{.5cm}}*{5}{M |}}\hhline{|-|-|-|-|-|}
  \coloredseed 1  & \cellcolor{white} 0  & \cellcolor{white} 0  & \cellcolor{white} 0  & \cellcolor{white} 0  \\\hhline{|-|-|-|-|-|}
  \cellcolor{white} 0  \\\hhline{|-|}
  \coloredseed 1  \\\hhline{|-|-|}
  \coloredseed 1  & \coloredseed 1  \\\hhline{|-|-|}
\end{tabular}}} \\\end{tabular}} ][.{\begin{tabular}{c}{{\bf \begin{tabular}{|@{\rule[-.15cm]{0pt}{.5cm}}*{5}{M |}}\hhline{|-|-|-|-|-|}
  \cellcolor{white} 0  & \cellcolor{white} 0  & \cellcolor{white} 0  & \cellcolor{white} 0  & \cellcolor{white} 0  \\\hhline{|-|-|-|-|-|}
  \cellcolor{white} 0  \\\hhline{|-|}
  \cellcolor{white} 0  \\\hhline{|-|}
  \coloredseed 1  \\\hhline{|-|-|}
  \coloredseed 1  & \coloredseed 1  \\\hhline{|-|-|}
\end{tabular}}} \\\end{tabular}} ]][.{\begin{tabular}{c}{{\bf \begin{tabular}{|@{\rule[-.15cm]{0pt}{.5cm}}*{4}{M |}}\hhline{|-|-|-|-|}
  \coloredseed 1  & \coloredseed 1  & \cellcolor{white} 0  & \coloredseed 1  \\\hhline{|-|-|-|-|}
  \coloredseed 1  & \coloredseed 1  \\\hhline{|-|-|}
\end{tabular}}} \\\end{tabular}} [.{\begin{tabular}{c}{{\bf \begin{tabular}{|@{\rule[-.15cm]{0pt}{.5cm}}*{4}{M |}}\hhline{|-|-|-|-|}
  \coloredseed 1  & \cellcolor{white} 0  & \coloredseed 1  & \coloredseed 1  \\\hhline{|-|-|-|-|}
  \coloredseed 1  & \coloredseed 1  & \coloredseed 1  \\\hhline{|-|-|-|}
\end{tabular}}} \\\end{tabular}} ][.{\begin{tabular}{c}{{\bf \begin{tabular}{|@{\rule[-.15cm]{0pt}{.5cm}}*{4}{M |}}\hhline{|-|-|-|-|}
  \cellcolor{white} 0  & \coloredseed 1  & \cellcolor{white} 0  & \coloredseed 1  \\\hhline{|-|-|-|-|}
  \cellcolor{white} 0  & \coloredseed 1  \\\hhline{|-|-|}
  \coloredseed 1  & \coloredseed 1  \\\hhline{|-|-|}
\end{tabular}}} \\\end{tabular}} ][.{\begin{tabular}{c}{{\bf \begin{tabular}{|@{\rule[-.15cm]{0pt}{.5cm}}*{4}{M |}}\hhline{|-|-|-|-|}
  \cellcolor{white} 0  & \cellcolor{white} 0  & \cellcolor{white} 0  & \cellcolor{white} 0  \\\hhline{|-|-|-|-|}
  \cellcolor{white} 0  & \cellcolor{white} 0  \\\hhline{|-|-|}
  \cellcolor{white} 0  \\\hhline{|-|}
  \coloredseed 1  \\\hhline{|-|-|}
  \coloredseed 1  & \coloredseed 1  \\\hhline{|-|-|}
\end{tabular}}} \\\end{tabular}} ]][.{\begin{tabular}{c}{{\bf \begin{tabular}{|@{\rule[-.15cm]{0pt}{.5cm}}*{4}{M |}}\hhline{|-|-|-|-|}
  \coloredseed 1  & \cellcolor{white} 0  & \cellcolor{white} 0  & \cellcolor{white} 0  \\\hhline{|-|-|-|-|}
  \coloredseed 1  \\\hhline{|-|-|}
  \coloredseed 1  & \coloredseed 1  \\\hhline{|-|-|}
\end{tabular}}} \\\end{tabular}} [.{\begin{tabular}{c}{{\bf \begin{tabular}{|@{\rule[-.15cm]{0pt}{.5cm}}*{4}{M |}}\hhline{|-|-|-|-|}
  \cellcolor{white} 0  & \cellcolor{white} 0  & \cellcolor{white} 0  & \coloredseed 1  \\\hhline{|-|-|-|-|}
  \coloredseed 1  \\\hhline{|-|-|-|}
  \coloredseed 1  & \coloredseed 1  & \coloredseed 1  \\\hhline{|-|-|-|}
\end{tabular}}} \\\end{tabular}} ]][.{\begin{tabular}{c}{{\bf \begin{tabular}{|@{\rule[-.15cm]{0pt}{.5cm}}*{4}{M |}}\hhline{|-|-|-|-|}
  \cellcolor{white} 0  & \cellcolor{white} 0  & \cellcolor{white} 0  & \cellcolor{white} 0  \\\hhline{|-|-|-|-|}
  \cellcolor{white} 0  \\\hhline{|-|}
  \coloredseed 1  \\\hhline{|-|-|}
  \coloredseed 1  & \coloredseed 1  \\\hhline{|-|-|}
\end{tabular}}} \\\end{tabular}} ]][.{\begin{tabular}{c}{{\bf \begin{tabular}{|@{\rule[-.15cm]{0pt}{.5cm}}*{3}{M |}}\hhline{|-|-|-|}
  \coloredseed 1  & \cellcolor{white} 0  & \coloredseed 1  \\\hhline{|-|-|-|}
  \coloredseed 1  & \coloredseed 1  \\\hhline{|-|-|}
\end{tabular}}} \\\end{tabular}} [.{\begin{tabular}{c}{{\bf \begin{tabular}{|@{\rule[-.15cm]{0pt}{.5cm}}*{3}{M |}}\hhline{|-|-|-|}
  \cellcolor{white} 0  & \coloredseed 1  & \coloredseed 1  \\\hhline{|-|-|-|}
  \coloredseed 1  & \coloredseed 1  & \coloredseed 1  \\\hhline{|-|-|-|}
\end{tabular}}} \\\end{tabular}} [.{\begin{tabular}{c}{{\bf \begin{tabular}{|@{\rule[-.15cm]{0pt}{.5cm}}*{3}{M |}}\hhline{|-|-|-|}
  \coloredseed 1  & \cellcolor{white} 0  & \coloredseed 1  \\\hhline{|-|-|-|}
  \coloredseed 1  & \cellcolor{white} 0  & \coloredseed 1  \\\hhline{|-|-|-|}
  \coloredseed 1  & \coloredseed 1  \\\hhline{|-|-|}
\end{tabular}}} \\\end{tabular}} ][.{\begin{tabular}{c}{{\bf \begin{tabular}{|@{\rule[-.15cm]{0pt}{.5cm}}*{3}{M |}}\hhline{|-|-|-|}
  \cellcolor{white} 0  & \cellcolor{white} 0  & \coloredseed 1  \\\hhline{|-|-|-|}
  \cellcolor{white} 0  & \cellcolor{white} 0  & \cellcolor{white} 0  \\\hhline{|-|-|-|}
  \coloredseed 1  \\\hhline{|-|-|}
  \coloredseed 1  & \coloredseed 1  \\\hhline{|-|-|}
\end{tabular}}} \\\end{tabular}} ]][.{\begin{tabular}{c}{{\bf \begin{tabular}{|@{\rule[-.15cm]{0pt}{.5cm}}*{3}{M |}}\hhline{|-|-|-|}
  \cellcolor{white} 0  & \cellcolor{white} 0  & \cellcolor{white} 0  \\\hhline{|-|-|-|}
  \cellcolor{white} 0  & \cellcolor{white} 0  \\\hhline{|-|-|}
  \coloredseed 1  \\\hhline{|-|-|}
  \coloredseed 1  & \coloredseed 1  \\\hhline{|-|-|}
\end{tabular}}} \\\end{tabular}} ]][.{\begin{tabular}{c}{{\bf \begin{tabular}{|@{\rule[-.15cm]{0pt}{.5cm}}*{3}{M |}}\hhline{|-|-|-|}
  \cellcolor{white} 0  & \cellcolor{white} 0  & \cellcolor{white} 0  \\\hhline{|-|-|-|}
  \coloredseed 1  \\\hhline{|-|-|}
  \coloredseed 1  & \coloredseed 1  \\\hhline{|-|-|}
\end{tabular}}} \\\end{tabular}} ]]\end{tikzpicture}}
    \end{center}
\caption{Seed tables of the children, grandchildren, and great grandchildren of  $\{0, 3, \dots\}$}\label{fig:exc}
\end{figure}

\begin{figure}[h]
   \begin{center}
     \input{Figure5.tex}
    \end{center}
  \caption{Seed tables of the children, grandchildren, and great grandchildren of  $\{0, 4, \dots\}$}\label{fig:exd}
\end{figure}

\begin{figure}[h]
\begin{center}
%
%
%
%

\newcolumntype{M}{>{\centering\arraybackslash}m{.5cm}}\setlength\tabcolsep{0pt}\setlength\arrayrulewidth{1pt}\providecommand\circledcolorednumb{}\renewcommand\circledcolorednumb[2]{\resizebox{0.025116\textwidth}{!}{\tikz[baseline=(char.center)]{\node[shape = circle,draw, inner sep = 2pt,fill=#1](char)    {\phantom{00}};\node[anchor=center] at (char.center) {\makebox(0,0){\large{#2}}};}}}
\robustify{\circledcolorednumb}
\providecommand\nongap{}\renewcommand\nongap[1]{\circledcolorednumb{yellow}{{\bf#1}}}
\providecommand\gap{}\renewcommand\gap[1]{\circledcolorednumb{white}{\phantom{#1}}}
\providecommand\generator{}\renewcommand\generator[1]{\circledcolorednumb{orange}{{\bf#1}}}
\providecommand\seed{}\renewcommand\seed[1]{\circledcolorednumb{blue!30}{{\bf#1}}}
\providecommand\nonseed{}\renewcommand\nonseed[1]{\circledcolorednumb{white}{{\bf#1}}}
\providecommand\dotscircles{}\renewcommand\dotscircles{\resizebox{0.025116\textwidth}{!}{\dots}}
\providecommand\gapingapset{}\renewcommand\gapingapset[1]{\circledcolorednumb{green!30}{#1}}
\providecommand\nongapingapset{}\renewcommand\nongapingapset[1]{\phantom{\gapingapset{#1}}}
\providecommand\coloredseed{}\renewcommand\coloredseed{\cellcolor{blue!30}}
\adjustbox{max width=\textwidth,max height=.9\textheight}{\begin{tikzpicture}[grow=down,sibling distance=10.000000mm]\tikzset{every tree node/.style={anchor=north}}\tikzset{level 1/.style={level distance=3.500000cm}}\tikzset{level 2/.style={level distance=4.375000cm}}\tikzset{level 3/.style={level distance=5.687500cm}}\tikzset{level 4/.style={level distance=7.000000cm}}\tikzset{level 5/.style={level distance=8.750000cm}}\tikzset{level 6/.style={level distance=8.925000cm}}\tikzset{level 7+/.style={level distance=9.450000cm}}\Tree[.{\begin{tabular}{c}{{\bf \begin{tabular}{|@{\rule[-.15cm]{0pt}{.5cm}}*{4}{M |}}\hhline{|-|-|-|-|}
  \coloredseed 1  & \cellcolor{white} 0  & \coloredseed 1  & \coloredseed 1  \\\hhline{|-|-|-|-|}
  \coloredseed 1  & \coloredseed 1  & \coloredseed 1  \\\hhline{|-|-|-|}
\end{tabular}}} \\\end{tabular}} [.{\begin{tabular}{c}{{\bf \begin{tabular}{|@{\rule[-.15cm]{0pt}{.5cm}}*{4}{M |}}\hhline{|-|-|-|-|}
  \cellcolor{white} 0  & \coloredseed 1  & \coloredseed 1  & \coloredseed 1  \\\hhline{|-|-|-|-|}
  \coloredseed 1  & \coloredseed 1  & \coloredseed 1  & \coloredseed 1  \\\hhline{|-|-|-|-|}
\end{tabular}}} \\\end{tabular}} [.{\begin{tabular}{c}{{\bf \begin{tabular}{|@{\rule[-.15cm]{0pt}{.5cm}}*{4}{M |}}\hhline{|-|-|-|-|}
  \coloredseed 1  & \coloredseed 1  & \cellcolor{white} 0  & \coloredseed 1  \\\hhline{|-|-|-|-|}
  \coloredseed 1  & \coloredseed 1  & \cellcolor{white} 0  & \coloredseed 1  \\\hhline{|-|-|-|-|}
  \coloredseed 1  & \coloredseed 1  \\\hhline{|-|-|}
\end{tabular}}} \\\end{tabular}} [.{\begin{tabular}{c}{{\bf \begin{tabular}{|@{\rule[-.15cm]{0pt}{.5cm}}*{4}{M |}}\hhline{|-|-|-|-|}
  \coloredseed 1  & \cellcolor{white} 0  & \coloredseed 1  & \coloredseed 1  \\\hhline{|-|-|-|-|}
  \coloredseed 1  & \cellcolor{white} 0  & \coloredseed 1  & \coloredseed 1  \\\hhline{|-|-|-|-|}
  \coloredseed 1  & \coloredseed 1  & \coloredseed 1  \\\hhline{|-|-|-|}
\end{tabular}}} \\\end{tabular}} ][.{\begin{tabular}{c}{{\bf \begin{tabular}{|@{\rule[-.15cm]{0pt}{.5cm}}*{4}{M |}}\hhline{|-|-|-|-|}
  \cellcolor{white} 0  & \coloredseed 1  & \cellcolor{white} 0  & \coloredseed 1  \\\hhline{|-|-|-|-|}
  \cellcolor{white} 0  & \coloredseed 1  & \cellcolor{white} 0  & \coloredseed 1  \\\hhline{|-|-|-|-|}
  \cellcolor{white} 0  & \coloredseed 1  \\\hhline{|-|-|}
  \coloredseed 1  & \coloredseed 1  \\\hhline{|-|-|}
\end{tabular}}} \\\end{tabular}} ][.{\begin{tabular}{c}{{\bf \begin{tabular}{|@{\rule[-.15cm]{0pt}{.5cm}}*{4}{M |}}\hhline{|-|-|-|-|}
  \cellcolor{white} 0  & \cellcolor{white} 0  & \cellcolor{white} 0  & \coloredseed 1  \\\hhline{|-|-|-|-|}
  \cellcolor{white} 0  & \cellcolor{white} 0  & \cellcolor{white} 0  & \cellcolor{white} 0  \\\hhline{|-|-|-|-|}
  \cellcolor{white} 0  & \cellcolor{white} 0  \\\hhline{|-|-|}
  \cellcolor{white} 0  \\\hhline{|-|}
  \coloredseed 1  \\\hhline{|-|-|}
  \coloredseed 1  & \coloredseed 1  \\\hhline{|-|-|}
\end{tabular}}} \\\end{tabular}} ]][.{\begin{tabular}{c}{{\bf \begin{tabular}{|@{\rule[-.15cm]{0pt}{.5cm}}*{4}{M |}}\hhline{|-|-|-|-|}
  \coloredseed 1  & \cellcolor{white} 0  & \cellcolor{white} 0  & \coloredseed 1  \\\hhline{|-|-|-|-|}
  \coloredseed 1  & \cellcolor{white} 0  & \cellcolor{white} 0  & \cellcolor{white} 0  \\\hhline{|-|-|-|-|}
  \coloredseed 1  \\\hhline{|-|-|}
  \coloredseed 1  & \coloredseed 1  \\\hhline{|-|-|}
\end{tabular}}} \\\end{tabular}} [.{\begin{tabular}{c}{{\bf \begin{tabular}{|@{\rule[-.15cm]{0pt}{.5cm}}*{4}{M |}}\hhline{|-|-|-|-|}
  \cellcolor{white} 0  & \cellcolor{white} 0  & \coloredseed 1  & \coloredseed 1  \\\hhline{|-|-|-|-|}
  \cellcolor{white} 0  & \cellcolor{white} 0  & \cellcolor{white} 0  & \coloredseed 1  \\\hhline{|-|-|-|-|}
  \coloredseed 1  \\\hhline{|-|-|-|}
  \coloredseed 1  & \coloredseed 1  & \coloredseed 1  \\\hhline{|-|-|-|}
\end{tabular}}} \\\end{tabular}} ][.{\begin{tabular}{c}{{\bf \begin{tabular}{|@{\rule[-.15cm]{0pt}{.5cm}}*{4}{M |}}\hhline{|-|-|-|-|}
  \cellcolor{white} 0  & \cellcolor{white} 0  & \cellcolor{white} 0  & \cellcolor{white} 0  \\\hhline{|-|-|-|-|}
  \cellcolor{white} 0  & \cellcolor{white} 0  & \cellcolor{white} 0  & \cellcolor{white} 0  \\\hhline{|-|-|-|-|}
  \cellcolor{white} 0  \\\hhline{|-|-|}
  \cellcolor{white} 0  & \cellcolor{white} 0  \\\hhline{|-|-|}
  \cellcolor{white} 0  \\\hhline{|-|}
  \coloredseed 1  \\\hhline{|-|-|}
  \coloredseed 1  & \coloredseed 1  \\\hhline{|-|-|}
\end{tabular}}} \\\end{tabular}} ]][.{\begin{tabular}{c}{{\bf \begin{tabular}{|@{\rule[-.15cm]{0pt}{.5cm}}*{4}{M |}}\hhline{|-|-|-|-|}
  \cellcolor{white} 0  & \cellcolor{white} 0  & \cellcolor{white} 0  & \coloredseed 1  \\\hhline{|-|-|-|-|}
  \cellcolor{white} 0  & \cellcolor{white} 0  & \cellcolor{white} 0  & \cellcolor{white} 0  \\\hhline{|-|-|-|-|}
  \cellcolor{white} 0  \\\hhline{|-|}
  \coloredseed 1  \\\hhline{|-|-|}
  \coloredseed 1  & \coloredseed 1  \\\hhline{|-|-|}
\end{tabular}}} \\\end{tabular}} [.{\begin{tabular}{c}{{\bf \begin{tabular}{|@{\rule[-.15cm]{0pt}{.5cm}}*{4}{M |}}\hhline{|-|-|-|-|}
  \cellcolor{white} 0  & \cellcolor{white} 0  & \cellcolor{white} 0  & \cellcolor{white} 0  \\\hhline{|-|-|-|-|}
  \cellcolor{white} 0  & \cellcolor{white} 0  & \cellcolor{white} 0  & \cellcolor{white} 0  \\\hhline{|-|-|-|-|}
  \cellcolor{white} 0  \\\hhline{|-|}
  \cellcolor{white} 0  \\\hhline{|-|-|}
  \cellcolor{white} 0  & \cellcolor{white} 0  \\\hhline{|-|-|}
  \cellcolor{white} 0  \\\hhline{|-|}
  \coloredseed 1  \\\hhline{|-|-|}
  \coloredseed 1  & \coloredseed 1  \\\hhline{|-|-|}
\end{tabular}}} \\\end{tabular}} ]]][.{\begin{tabular}{c}{{\bf \begin{tabular}{|@{\rule[-.15cm]{0pt}{.5cm}}*{4}{M |}}\hhline{|-|-|-|-|}
  \coloredseed 1  & \cellcolor{white} 0  & \cellcolor{white} 0  & \coloredseed 1  \\\hhline{|-|-|-|-|}
  \cellcolor{white} 0  & \cellcolor{white} 0  & \cellcolor{white} 0  \\\hhline{|-|-|-|}
  \coloredseed 1  \\\hhline{|-|-|}
  \coloredseed 1  & \coloredseed 1  \\\hhline{|-|-|}
\end{tabular}}} \\\end{tabular}} [.{\begin{tabular}{c}{{\bf \begin{tabular}{|@{\rule[-.15cm]{0pt}{.5cm}}*{4}{M |}}\hhline{|-|-|-|-|}
  \cellcolor{white} 0  & \cellcolor{white} 0  & \coloredseed 1  & \cellcolor{white} 0  \\\hhline{|-|-|-|-|}
  \cellcolor{white} 0  & \cellcolor{white} 0  & \coloredseed 1  \\\hhline{|-|-|-|}
  \coloredseed 1  \\\hhline{|-|-|-|}
  \coloredseed 1  & \coloredseed 1  & \coloredseed 1  \\\hhline{|-|-|-|}
\end{tabular}}} \\\end{tabular}} [.{\begin{tabular}{c}{{\bf \begin{tabular}{|@{\rule[-.15cm]{0pt}{.5cm}}*{4}{M |}}\hhline{|-|-|-|-|}
  \cellcolor{white} 0  & \cellcolor{white} 0  & \cellcolor{white} 0  & \coloredseed 1  \\\hhline{|-|-|-|-|}
  \cellcolor{white} 0  & \cellcolor{white} 0  & \cellcolor{white} 0  \\\hhline{|-|-|-|}
  \coloredseed 1  \\\hhline{|-|-|-|}
  \cellcolor{white} 0  & \cellcolor{white} 0  & \cellcolor{white} 0  \\\hhline{|-|-|-|}
  \coloredseed 1  \\\hhline{|-|-|}
  \coloredseed 1  & \coloredseed 1  \\\hhline{|-|-|}
\end{tabular}}} \\\end{tabular}} ]][.{\begin{tabular}{c}{{\bf \begin{tabular}{|@{\rule[-.15cm]{0pt}{.5cm}}*{4}{M |}}\hhline{|-|-|-|-|}
  \cellcolor{white} 0  & \cellcolor{white} 0  & \cellcolor{white} 0  & \cellcolor{white} 0  \\\hhline{|-|-|-|-|}
  \cellcolor{white} 0  & \cellcolor{white} 0  & \cellcolor{white} 0  \\\hhline{|-|-|-|}
  \cellcolor{white} 0  \\\hhline{|-|-|}
  \cellcolor{white} 0  & \cellcolor{white} 0  \\\hhline{|-|-|}
  \cellcolor{white} 0  \\\hhline{|-|}
  \coloredseed 1  \\\hhline{|-|-|}
  \coloredseed 1  & \coloredseed 1  \\\hhline{|-|-|}
\end{tabular}}} \\\end{tabular}} ]][.{\begin{tabular}{c}{{\bf \begin{tabular}{|@{\rule[-.15cm]{0pt}{.5cm}}*{4}{M |}}\hhline{|-|-|-|-|}
  \cellcolor{white} 0  & \cellcolor{white} 0  & \cellcolor{white} 0  & \cellcolor{white} 0  \\\hhline{|-|-|-|-|}
  \cellcolor{white} 0  & \cellcolor{white} 0  & \cellcolor{white} 0  \\\hhline{|-|-|-|}
  \cellcolor{white} 0  \\\hhline{|-|}
  \coloredseed 1  \\\hhline{|-|-|}
  \coloredseed 1  & \coloredseed 1  \\\hhline{|-|-|}
\end{tabular}}} \\\end{tabular}} ]]\end{tikzpicture}}
  \end{center}
  \caption{Seed tables of the children, grandchildren, and great grandchildren of  $\{0, 4, 7, \dots\}$}\label{fig:exe}
\end{figure}


\begin{figure}[h]
  \input{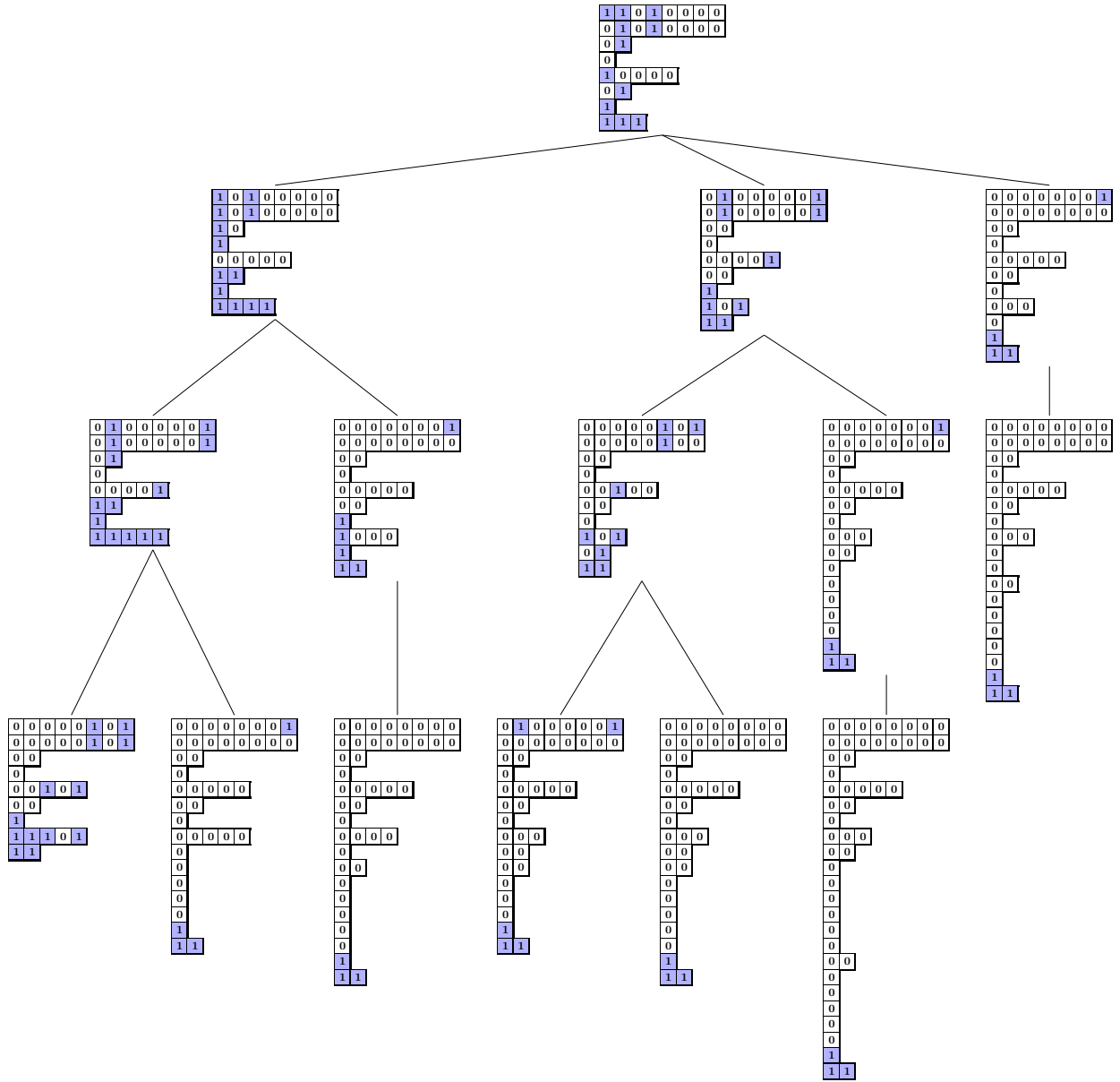}
  \caption{Seed tables of the children, grandchildren, and great grandchildren of  $\{0, 8, 16, 18, 19, 24, 26, 27, 30, \dots\}$}\label{fig:exf}
\end{figure}

\newpage
\section{Implementation of the improved seeds algorithm}
\lstset{language=C}

\enlargethispage{2cm}

{\tiny
\begin{lstlisting}[numbers=left,firstnumber=1,firstline=18,lastline=500]
#include <stdio.h>
#include <stdlib.h>
#include <time.h>
typedef unsigned __int128 bbint;
int gamma;
long long int Sdesc(const bbint G, bbint S, const int c, const int m, const int u, const int v, int r, int gd) {
  int newc, s = 1, sold = 0;
  bbint A;
  long long int ng = 0;
  if (gd) {
    --gd;
    A = G;
    if (S & (bbint)1) {
      if (S & ((bbint)1 << m))
        ng += Sdesc(G | ((bbint)1 << (c - 1)), (S >> 1) | ((bbint)3 << (c - 1)), c + 1, m, u, v, r--, gd);
      else
        ng += Sdesc(G | ((bbint)1 << (c - 1)), (S >> 1) | ((bbint)3 << (c - 1)), c + 1, m, u, v, --r, gd);
    }
    while (r) {
      if (S & (bbint)1 << s) {
        newc = c + s + 1;
        for (int i = sold; i < s; i++) {
          A <<= 1;
          S &= A;
        }
        if (S & ((bbint)1 << (m + s)))
          ng += Sdesc(G | ((bbint)1 << (newc - 2)), S >> (s + 1) | ((bbint)7 << (newc - 3)), newc, m, u, v, r--, gd);
        else
          ng += Sdesc(G | ((bbint)1 << (newc - 2)), S >> (s + 1) | ((bbint)7 << (newc - 3)), newc, m, u, v, --r, gd);
        sold = s;
      }
      s++;
    }
  } else {
    A = S & (S >> m) & (((bbint)1 << u) - 1);
    s = 0;
    while (A) {
      s++;
      A &= A - 1;
    }
    ng += s * (r - 1);
    ng += r * (r - 1) * (r - 2) / 6;
    if (c > m + u) {
      A = (((bbint)1 << v) - 1) & S & (S >> u) & (S >> (m + u));
      while (A) {
        ng++;
        A &= A - 1;
      }
    }
  }
  return ng;
}
long long int lowrank(const int m, const int u, const int gamma) {
  int gd, min, mu, c;
  bbint G, S;
  long long int ng = 0;
  mu = m + u;
  if (u == m) {
    gd = gamma - mu;
    if (m < gd)
      min = m;
    else
      min = gd;
    c = mu;
    for (int v = 1; v < min; v++) {
      gd--;
      G = ((bbint)1 << c) - 1;
      c++;
      G -= ((bbint)1 << (m - 1));
      G -= ((bbint)1 << (mu - 1));
      S = ((bbint)1 << c) - 1;
      S -= ((bbint)1 << (m - v));
      S -= ((bbint)1 << (mu - v));
      ng += Sdesc(G, S - 1, c, m, u, v, m - 2, gd);
    }
    gd--;
    G = ((bbint)1 << c) - 1;
    c++;
    G -= ((bbint)1 << (m - 1));
    G -= ((bbint)1 << (mu - 1));
    S = ((bbint)1 << c) - 1;
    S -= ((bbint)1 << (m - min));
    S -= ((bbint)1 << (mu - min));
    ng += Sdesc(G, S, c, m, u, min, m - 1, gd);
  } else {
    if (m - u < gamma - mu)
      min = m - u;
    else
      min = gamma - mu;
    c = mu;
    gd = gamma - m - u;
    for (int v = 1; v < min; v++) {
      gd--;
      G = ((bbint)1 << c) - 1;
      c++;
      G -= ((bbint)1 << (m - 1));
      G -= ((bbint)1 << (mu - 1));
      S = ((bbint)1 << c) - 1;
      S -= ((bbint)1 << (m - u - v));
      S -= ((bbint)1 << (m - v));
      S -= ((bbint)1 << (mu - v));
      if (u < v)
        ng += Sdesc(G, S - 1, c, m, u, v, m - 4, gd);
      else
        ng += Sdesc(G, S - 1, c, m, u, v, m - 3, gd);
    }
    gd--;
    G = ((bbint)1 << c) - 1;
    c++;
    G -= ((bbint)1 << (m - 1));
    G -= ((bbint)1 << (mu - 1));
    S = ((bbint)1 << c) - 1;
    S -= ((bbint)1 << (m - u - min));
    S -= ((bbint)1 << (m - min));
    S -= ((bbint)1 << (mu - min));
    if (u < min)
      ng += Sdesc(G, S, c, m, u, min, m - 3, gd);
    else
      ng += Sdesc(G, S, c, m, u, min, m - 2, gd);
  }
  return ng;
}
int main(int numvars, char **vars) {
  long long int ng = 0;
  bbint G, S;
  int x, min;
  time_t seconds, secondsafter;
  gamma = (int)atoi(vars[1]);
  seconds = time(NULL);
  ng = (gamma - 4) * (gamma - 5) * (gamma - 6) / 6 + (gamma - 2) * (gamma - 3) + 8; // case m=gamma-2 and the hyperelliptic;
  for (int m = 3; m < gamma - 2; m++) {
    x = gamma - m - 1;
    if (m - 1 < x)
      min = m - 1;
    else
      min = x;
    ng += Sdesc(((bbint)1 << (m - 1)) - 1, ((bbint)1 << m) - 8, m, m, 1, 1, m - 3, x - 1);
    for (int v = 2; v < min; v++) {
      G = ((bbint)1 << (m + v)) - 1;
      G -= ((bbint)3 << (m - 1));
      S = ((bbint)1 << (m + v + 1)) - 1;
      S -= ((bbint)7 << (m - v - 1));
      ng += Sdesc(G, S - 1, m + v + 1, m, 1, v, m - 4, x - v);
    }
    G = ((bbint)1 << (m + min)) - 1;
    G -= ((bbint)3 << (m - 1));
    S = ((bbint)1 << (m + min + 1)) - 1;
    S -= ((bbint)7 << (m - min - 1));
    ng += Sdesc(G, S, m + min + 1, m, 1, min, m - 3, x - min);
    if (x > m) {
      for (int u = 2; u <= m; u++)
        ng += lowrank(m, u, gamma);
    } else {
      for (int u = 2; u < x; u++) {
        ng += lowrank(m, u, gamma);
      }
      ng += (m - 1) * (m - 2) * (m - 3) / 6;
      if (m < 2 * x)
        ng += (x - 1) * (m - 2);
      else
        ng += x * (m - 2);
      if (x == m || 2 * x != m)
        ng++;
      if (x == m - 1)
        ng++;
      if (x < m - 1 && 2 * x != m - 1)
        ng += 2;
    }
  }
  secondsafter = time(NULL);
  printf("n%d=%lld (without parallelization) time taken %d\n", (int)gamma, ng, (int)(secondsafter - seconds));
  return 0;
}
\end{lstlisting}
}


\end{document}